\documentclass[11pt,letterpaper]{article}
\usepackage[dvips,linkcolor=blue,urlcolor=blue,citecolor=blue,linktocpage,bookmarks,pdfstartview=FitBH, breaklinks = true]{hyperref}
\usepackage{fancyhdr}
\usepackage[all]{hypcap}
\usepackage{latexsym}
\usepackage{verbatim}
\usepackage{setspace}
\usepackage[centertags]{amsmath}
\usepackage{amsfonts}
\usepackage{amssymb}
\usepackage{amsthm}
\usepackage{mathrsfs}
\usepackage{newlfont}
\usepackage{natbib}
\usepackage{graphicx}
\usepackage{lscape}
\usepackage{slashbox}
\usepackage{bm}
\usepackage{appendix}
\usepackage{rotating}
\usepackage{url}
\usepackage{float}
\usepackage{enumerate}
\usepackage{datetime}
\usepackage{fancyvrb} 
\usepackage{authblk} 
\usepackage{lineno} 

\numberwithin{equation}{section}

\theoremstyle{plain}

\theoremstyle{definition}

\theoremstyle{remark}

\numberwithin{equation}{section}

\def\0{\mathbf{\dobold 0}}

\def\1{\mathbf{\dobold 1}}

\def\hbbetaS+{\bm{\hat\beta}^{\textrm{S+}}}

\def\D2{\Delta^2} 

\title{Improved LASSO}

\author[1]{A. K. Md. Ehsanes Saleh}
\author[2]{Enayetur Raheem}
\affil[1]{Carleton University, Ottawa, ON, Canada}
\affil[2]{University of Northern Colorado, Greeley, CO, USA}

\date{}

\begin{document}

\maketitle


\centerline {Version 1.4 updated  \today~ \currenttime}

\vspace {.5in}

\centerline{\bf Abstract}

\bigskip
We propose an improved LASSO estimation technique based on Stein-rule. We shrink classical LASSO estimator using preliminary test, shrinkage, and positive-rule shrinkage principle. Simulation results have been carried out for various configurations of correlation coefficients ($r$), size of the parameter vector ($\beta$), error variance ($\sigma^2$) and number of non-zero coefficients ($k$) in the model parameter vector. Several real data examples have been used to demonstrate the practical usefulness of the proposed estimators. Our study shows that the risk ordering  given by LSE $>$ LASSO $>$ Stein-type LASSO $>$ Stein-type positive rule LASSO, remains the same uniformly in the divergence parameter $\Delta^2$ as in the traditional case.



\section{Introduction}\label{sec:intro}

The estimation of parameters of a model with ``uncertain prior information'' on parameters of interest began with \cite{bancroft:1944} in the classical front. But a breakthrough came when \cite{stein:1956} and \cite{james:stein:1961} proved that the sample mean in a multivariate normal model is not admissible under a quadratic loss, for dimension more than two. This very result gave birth to a class of shrinkage estimators of various form and setup. A partial document on preliminary test and Stein-type estimators are given by \cite{judge:bock:1978}. The Stein-type estimators have been reformulated and expanded by \citet[Ch 4.4.3]{saleh:2006} which includes asymptotic and nonparametric methods. Due to the immense impact of Stein's results, scores of technical papers appeared in the literature covering various areas of applications. Here is one with the popular LASSO estimator.

This paper is devoted to the study of the performance characteristics of several improved estimators of LASSO based on preliminary test and Stein's principle, and the comparison of the LASSO estimator with the least square estimator (LSE), improved preliminary test estimator (IPT), shrinkage preliminary test estimator (SPTE), Stein-type shrinkage LASSO estimator (SLE), and  Stein-type positive-rule shrinkage LASSO estimator (PSLE). An important characteristic of LASSO-type estimators is that they provide simultaneous estimation and selection of coefficients in  linear models and can be applied when dimension of the parameter space exceeds the dimension of the sample space. Our conclusions are made based on simulated mean-squared errors (MSE) and relative efficiency tables and graphs. It is shown that the modified LASSO carries on with the same dominance characteristics as the usual preliminary test and Stein-type
estimators \citep{saleh:2006}.

Our contribution in this paper is unique. We have proposed a set of LASSO-based shrinkage estimators that perform superior to the classical LASSO estimator. We studied the theoretical properties of the estimators in terms of asymptotic mean squared errors (AMSE). Analytical expressions for the asymptotic risk functions of the proposed estimators have been provided. We carried out Monte Carlo simulation experiments to study the risk-behavior of the proposed estimators and their comparisons with the LASSO estimator. Application of the proposed estimators have been demonstrated using three real life data examples where the proposed estimators performed superiorly to the classical LASSO estimator. 

The organization of the paper is as follows. Section 2 contains the basic informations about the LASSO and preliminary test and Stein-type estimators. Proposed improved estimators of LASSO are presented here. Risk properties and risk-comparisons of various estimators are presented in this section. In section 3 we outline our simulation setup and discuss simulation results. Details of the analysis of relative efficiencies of estimators are presented in this section. In section 4, applications of the proposed estimators have been demonstrated using three real life data sets. Finally, conclusions are provided in section 5.

\section{Linear model and estimators}\label{sec:two}
Consider the linear multiple regression model 
\begin{equation}\label{eq2.1}
Y=X\beta + e
\end{equation}
where $Y=(y_1, y_2, \ldots, y_n)',$ $\beta=(\beta_1, \beta_2, \ldots, \beta_p)',$ $X$ is the $n \times p$  design matrix and 
$e=(e_1, e_2, \ldots, e_n)'$ iid error $n$-vector. Further we assume that $E(e)=0$ and $E(e'e)=\sigma^2 I_n$.

It is well-known that the LSE of $\beta$ is given by 
\begin{equation}\label{eq:beta:ols}
\tilde{\beta}_n = (X'X)^{-1}X'Y
\end{equation}
which we use for obtaining preliminary test and shrinkage estimators.  LSE $\tilde{\beta}_n$ is the ``best linear unbiased estimator (BLUE)'' of $\beta$. The solution (\ref{eq:beta:ols}) depends on the non-singularity of the matrix $X'X$. If it is singular, the uniqueness of $\hat{\beta}_n$ is lost and we end up with multiple solutions with varying variances and some of them may be very large. To avoid these problems, a class of penalty estimators evolved in the class of restricted estimators. A simple example of ``restricted estimator'' when we want an estimator of $\beta$ which belongs to the subspace defined by $H\beta=h$, where $H$ is a $q\times p$ matrix and $h$ is a $q$-vector or real numbers. 

Next, we consider the classical penalty estimator, called the LASSO (least absolute shrinkage and selection operator) belonging to a class of restricted estimators.

The estimators of $\beta$ is obtained by minimizing the LS criterion subject to $H\beta=h$. Explicitly we may write 
\[
\mathop{\min}_{\beta \in R^p} (Y-X\beta)'(Y-X\beta) + \lambda(H\beta-h).
\]
The solution for this problem is the ``restricted estimator'' and the tuning parameter $\lambda$ can be explicitly obtained giving
\begin{equation}\label{eq2.4}
\hat{\beta}_n = \tilde{\beta}_n -C^{-1}H'(HC^{-1}H')^{-1}(H\tilde{\beta}_n-h), \quad C=X'X,
\end{equation}
where $\tilde{\beta}_n = (X'X)^{-1}X'Y$, the LSE.

For this, we consider $\sum_{j=1}^{p}|\beta_j| < t$ as the restriction and minimize
\begin{equation}\label{eq:2.6}
\mathop{\min}_{\beta \in R^p} (Y-X\beta)'(Y-X\beta) + \lambda\sum_{j=1}^p |\beta_j|
\end{equation}
yielding the solution as LASSO introduced by \cite{tibshirani:1996}, and is given by
\begin{equation}\label{eq:2.7}
\hat{\beta}^L_n = \left\{
\begin{array}{l}
\sum_{i=1}^{n} \left(\sum_{j=1}^{p} x_{ij}x_{ik}\beta_j - x_{ik}y_i\right)+\frac{\lambda}{2}\mbox{sgn}(\beta_k)=0 \\
\mbox{ and } \beta_k=0
\end{array} 
\right.
\end{equation}
When $\frac{1}{n}X'X \rightarrow I_p$, the solution may be written as 
\begin{align}\label{eq:2.8}
\hat{\beta}_n^L &= \left(\hat{\beta}_{1n}^L,\ldots, \hat{\beta}_{pn}^L\right)\\
\nonumber
\hat{\beta}_{jn}^{L} &= \mbox{ sgn}(\tilde{\beta}_{jn})\left(|\tilde{\beta}_{jn}| -\frac{\lambda}{2}\right)^{+}, j=1, \ldots, p
\end{align} where $|\tilde{\beta}_n|=\left(|\tilde{\beta}_{1n}|, \ldots, |\tilde{\beta}_{pn}|\right)'.$

Actually \cite{frank:friedman:1993} defined the class of generalized version of LASSO, namely, the ``bridge estimator'' as 
\begin{eqnarray}\label{eq:2.9}
\mathop{\min}_{\beta\in R^p}(Y-X\beta)'(Y-X\beta) + \lambda\sum_{j=1}^{p}|\beta_j|^\gamma.
\end{eqnarray}
If $\gamma=1$ the solution reduces to the LASSO estimator.

LASSO proposed by \cite{tibshirani:1996} simultaneously estimates and makes selection of variables with appropriate interpretation and its viral popularity in applications. For computational solution and methodology see \cite{tibshirani:1996} and \cite{efron:etal:2004}. Later \cite{efron:etal:2004} proposed Least Angle Regression (LAR) which is a stepwise regression, and \cite{friedman:hastie:tibshirani:2010} developed an efficient algorithm for the estimation of a GLM with convex penalty. During the course of development of penalty estimators, \cite{fan:li:2001} defined good penalty functions as the one which yield (i) nearly unbiased estimator when true parameter is large to avoid unnecessary modeling bias, (ii) an estimator which is a threshold rule that sets small estimated coefficients to zero to reduce model complexity, and (iii) the resulting estimator to be continuous in the data to avoid instability in the model prediction. In this paper, we present an improved version of LASSO.

\subsection{PTE and Stein-type estimators}
For the linear multiple regression model, $Y=X\beta+e$, if we suspect the full hypothesis to be $\beta=0$ (null-vector), then the restricted estimator (RE) $\hat{\beta}_n=0$ and the test for $\beta=0$, vs $\beta \ne 0$ may be based on the statistic
\begin{equation}\label{eq:test:statistic}
\mathcal{L}_n = \frac{\tilde{\beta}_n'C\tilde{\beta}_n}{s_e^2},
\end{equation}
where 
\begin{equation}\label{eq:sigma2}
s_n^2 = (n-p)^{-1}(Y-X\tilde{\beta}_n)'(Y-X\tilde{\beta}_n).
\end{equation}
Under the conditions
\begin{enumerate}[(1)]
\item $\frac{1}{n}X'X \rightarrow C \mbox{ as } n\rightarrow \infty \, (C \mbox{ positive definite })$, and
\item $\mathop{\max}_{1 \leq j \le n} \bm{x}_j'(\frac{1}{n}X'X)^{-1}\bm{x}_j \rightarrow 0 \mbox{ as } n \rightarrow \infty$
\end{enumerate}
where $\bm{x}_j$ is the $j$th row of $X$,  $\mathcal{L}_n \rightarrow\chi^2_p$ --central chi-square variable with $p$ degrees of freedom (df).

Let $\chi^2_p(\alpha)$ be an upper $\alpha$-level critical value from this null distribution; then we may define the preliminary test estimator (PTE) of $\beta$ as
\begin{equation}\label{eq:2.14}
\hat{\beta}_n^{PT} = \tilde{\beta}_n -\tilde{\beta}_n \, I(\chi^2 < \chi^2_p(\alpha)).
\end{equation}
The PTE is a discrete variable. As a result some optimality properties when we consider assessing its MSE comparison is lost. We may define a continuous version of PTE as the James-Stein-type estimator (JSE) given by
\begin{equation}\label{eq:2.15}
\hat{\beta}_n^{S}=\tilde{\beta}_n -(p-2)\tilde{\beta}_n \mathcal{L}_n^{-1}.
\end{equation}
Note that we have replaced $I(\chi^2 < \chi^2_p(\alpha))$ by $(p-2)\mathcal{L}_n^{-1}$ in the definition of PTE. However, $\hat{\beta}_n^{S}$ has an inherent problem of changing its sign due to the factor $(1-(p-2)\mathcal{L}_n^{-1})$ which may be larger than 1 in absolute value. If that happens, from applied point of view, its interpretation becomes blurred. Thus, we define another estimator, namely, the positive-rule Stein-type estimator (PRSE) as
\begin{equation}\label{eq:2.16}
\hat{\beta}_n^{S+} = \hat{\beta}_n^{S} \, I(\mathcal{L}_n > (p-2)).
\end{equation}
Next, we define an improved preliminary test (IPT) estimator defined by 
\begin{equation}\label{eq:2.17}
\hat{\beta}_n^{IPT} = \hat{\beta}_n^{PT}\, \left(1-(p-2)\mathcal{L}_n^{-1} \right).
\end{equation}
Thus, to set the stage, we have defined six estimators, namely, LSE, RE, PTE, JSE, PRSE, and IPT here, and one penalty estimator LASSO. Next, we use the definitions above to propose new shrinkage-type LASSO estimators.

\subsection{Proposed improved estimators of LASSO}

Let us redefine the LASSO estimator (LE) \citep{tibshirani:1996} as 
\begin{equation}\label{eq:lasso}
\hat{\beta}_n^{LE} = \underset{\beta}{\mbox{argmin}}\, (y-X\beta)'(y-X\beta) + \lambda \sum_{j=1}^{p}|\beta_j|.
\end{equation}
From now on, we consider $\hat{\beta}_n^{LE}$ as our unrestricted estimator (UE). Then, similar to the definition of PTE, we define the preliminary test LASSO estimator (PTLE) as
\begin{equation}\label{eq:PTLE}
\hat{\beta}_n^{PTLE} = \hat{\beta}_n^{LE} \, I(\mathcal{L}_n \ge c_\alpha), \quad 0 \leq \alpha \leq 1
\end{equation}
The Stein-type shrinkage LASSO estimator (SLE) based on $\hat{\beta}_n^{LE}$ may be defined as
\begin{equation}\label{eq:SLE}
\hat{\beta}_n^{SLE} = \hat{\beta}_n^{LE} \left(1-(p-2) \mathcal{L}_n^{-1}\right),
\end{equation}
where $\mathcal{L}_n$ was defined in (\ref{eq:test:statistic}).
Now we define the Stein-type positive rule shrinkage LASSO estimator (PSLE) using $\hat{\beta}_n^{LE}$ and $\mathcal{L}_n$ as
\begin{equation}\label{eq:PSLE}
\hat{\beta}_n^{PSLE} = \hat{\beta}_n^{SLE} \, I(\mathcal{L}_n > (p-2)).
\end{equation}
We note that $\mathcal{L}_n$ follows a non-central $\chi^2$ distribution with noncentrality parameter $\Delta^2=\frac{\beta'C\beta}{\sigma^2}$ under local alternatives. 

\subsection{Orthogonal case, $r=0$}

In this section, we consider the model (\ref{eq2.1}) and assume that the design matrix is centered and $\frac{1}{n}(X'X)=I_p$. Under this condition, the LASSO estimator is given by 
\begin{equation}\label{eq:lasso:orthogonal:r0}
\hat{\beta}_n^{LE}(\lambda) = (\hat{\beta}_{1n}^{LE}, \ldots, \hat{\beta}_{pn}^{LE})',
\end{equation}
where $\hat{\beta}_{jn}^{LE}=\mbox{sgn }(\tilde{\beta}_{jn}) (|\tilde{\beta}_{jn}| - \lambda) I(|\beta_{jn| > \lambda})$ for $j = 1, \ldots, p$.

It is known from \cite{donoho:1994:ideal} as $n \rightarrow \infty$, the quadratic risk bound is given by 
\[
m = \frac{\sigma^2}{p}(1 + 2n p) (1 + \Delta_{\mbox{min}}^2)
\]
where $\Delta^2_{\mbox{min}} = \sum_{j=1}^{p}(\frac{\beta_j^2}{\sigma^2}, 1)$ under the local alternative 
\[
\beta_{jn} = n^{-1/2} \bm{\delta}, \quad \bm{\delta}=(\delta_1, \ldots, \delta_p)' \neq 0 \quad \mbox{for } \lambda=\sigma\sqrt{2np}.
\]
For the sparse solution, one has to use $\Delta_{\min}=k$ as such $m$ equals $m^*=\frac{1}{p}(1+2\ln p)(1+k)$ since we have $k (<p)$ coefficients satisfy $\delta^2_j>\sigma^2$ and remaining equal to zero.

Now, we consider the PTE, SE and PRSE as defined in section (2.1) where we take $\mathcal{L}_n=\frac{\tilde{\beta}_n'\tilde{\beta}_n}{s^2_n}$. Thus, one can find that the asymptotic risk-bound of the PTE, SE and PRSE are then given by
\begin{align*}
R_1(\hat{\beta}_n^{LE}) &= \frac{\sigma^2}{p}(1+2 \ln p) (1+k)=\sigma^2 m^*\\
R_2(\hat{\beta}_n)&=R_2(\hat{\beta}_n^{LE}) = \frac{\delta'\delta}{\sigma^2} = \Delta^2  \mbox{ Note that restricted estimator is zero.} \\
R_3(\hat{\beta}_n^{PTLE}) &= \sigma^2 \left[m^* (1-\mathcal{H}_{k+2}(c_\alpha; \Delta^2)) + \Delta^2\left\{2\mathcal{H}_{k+2} (c_\alpha;\Delta^2) - \mathcal{H}_{k+4}(c_\alpha;\Delta^2)\right\}\right] \\
R_4(\hat{\beta}_n^{SLE}) &= \sigma^2 \left[m^* - (p-2)m^* \left\{2E[\chi^{-2}_{k+2}(\Delta^2)]-(k-2)E[\chi^{-4}_{k+2}(\Delta^2)]\right\}\right. \\
&  \quad \left.  +(p^2-4)\Delta^2 E[\chi^{-4}_{k+4}(\Delta^2)]\right]\\
R_5(\hat{\beta}_n^{PSLE}) &= R_4(\hat{\beta}_n^{SLE}) -\sigma^2 m^* E[(1-(k-2)\chi^{-2}_{k+2}(\Delta^2))^2 I(\chi^2_{k+2}(\Delta^2)< k-2)] \\
& \quad +\sigma^2\Delta^2\left\{2E[(1-(k-2)\chi^{-2}_{k+4}(\Delta^2))^2 I(\chi^2_{k+4}(\Delta^2)<k-2)]\right. \\
& \quad \left. -E[(1-(k-2)\chi^{-2}_{k+4}(\Delta^2))^2 I (\chi^2_{k+4}(\Delta^2)<k-2)\right\}
\end{align*}
Here, $\mathcal{H}_{p+2\nu}(c_\alpha;\Delta^2)$ is the cdf of a noncentral chi-square distribution with $p+2\nu$ df and noncentrality parameter $\Delta^2,$ and 
\[
E[\chi^{-2r}_{p+2\nu} (\Delta^2)] =\int_0^\infty x^{-2r}d\mathcal{H}_{p+2\nu}(x; \Delta^2).
\]

\subsection{Analysis of asymptotic MSE of the estimators}
First, we note that $R_1(\tilde{\beta}_n) \ge R_1(\hat{\beta}_n^{LE})$ uniformly in $\Delta^2$. Next, we compare $\hat{\beta}_n^{LE}$ and $\hat{\beta}_n^{PTLE}$ by taking the risk-difference
\begin{align*}
R_1(\hat{\beta}_n^{LE}) - R_3(\hat{\beta}_n^{PTLE}) &= m^* \mathcal{H}_{k+2}(c_\alpha;\Delta^2)-\Delta^2\left\{2\mathcal{H}_{k+2}(c_\alpha;\Delta^2)-\mathcal{H}_{k+4}(c_\alpha;\Delta^2)\right\} \gtreqqless 0
\end{align*}
whenever 
\[
\Delta^2 \lesseqqgtr \frac{m^* \mathcal{H}_{k+2}(c_\alpha;\Delta^2)}{\left\{2\mathcal{H}_{k+2}(c_\alpha;\Delta^2)-\mathcal{H}_{k+4}(c_\alpha;\Delta^2)\right\}}.
\]
Note 
\begin{align*}
& m^* \mathcal{H}_{k+2} > \Delta^2(2 \mathcal{H}_{k+2}-\mathcal{H}_{k+4}) \\
\mbox{or, } & \Delta^2 < \frac{m^* \mathcal{H}_{k+2}}{2 \mathcal{H}_{k+2}-\mathcal{H}_{k+4}}
\end{align*}
Accordingly, PTLE is better than LE in the range of $\Delta^2$. Otherwise, for 
\[
\Delta^2 > \frac{m^* \mathcal{H}_{k+2}}{2 \mathcal{H}_{k+2}-\mathcal{H}_{k+4}}
\]
LE is better than PTLE. Next, we note that for all $\Delta^2$, the risk-differences 
\[
(i) \quad R_1(\hat{\beta}_n^{LE})-R_4(\hat{\beta}_n^{SLE}) \ge 0\quad \mbox{ for all } \Delta^2
\]
and risk-difference
\[
(ii) \quad R_4(\hat{\beta}_n^{SLE})-R_5(\hat{\beta}_n{PSLE}) \ge 0 \quad \mbox{for all } \Delta^2 \in (0, \infty).
\]
Hence, 
\[
R_5(\hat{\beta}_n^{PSLE}) \leq R_4(\hat{\beta}_n^{SLE}) \leq R_1(\hat{\beta}_n^{LE}) \quad \forall \Delta^2.
\]

\section{Simulation for orthogonal case} \label{sec:simulation}

For the non-orthogonal case, we conduct Monte Carlo simulation experiments to study the performance of the proposed positive-rule shrinkage-LASSO estimator (PSLE) along with preliminary test LASSO estimator (PTLE), and shrinkage LASSO estimator (SLE). In particular, we study relative efficiencies of the proposed estimators compared to the LASSO estimator (LE)  by \cite{tibshirani:1996}. In the simulation studies, mean squared errors (MSE) were computed for each of the proposed estimators, and their relative efficiencies were calculated by taking the ratio of MSE of the proposed estimators to the MSE of LE. \cite{raheem:ahmed:doksum:2012} have conducted similar studies where a sub-hypothesis was tested and relative efficiencies of various shrinkage and penalty estimators were studied in a partially linear regression setup. In this study, we are concerned with full model hypothesis $H_0:\beta=0$ against the alternative $H_a: \beta \ne 0$. \cite{saleh:raheem:2015} have studied the performance of various shrinkage and penalty estimators under a full model hypothesis using the same setup. Next, we discuss the simulation setup. 

First, We generate the design matrix $\bm{X}$ from a multivariate normal distribution with mean vector $\bm{\mu} =\bm{0}$ and covariance matrix $\bm{\Sigma}$. The off-diagonal elements of the covariance matrix are considered to be equal to $r$ with $r=0, 0.5, 0.9$. We consider sample size to be $n=100$, and the number of parameters, $p$, equal to 10, 20, and 50. In our setup, $\beta$ is a $p$-vector and a function of $\Delta^2$. In the simulation, the $\beta$ vector is defined such that a $\Delta^2=0$ indicates a data set being generated under null hypothesis, whereas $\Delta^2 >0$ indicates a data set generated under alternative hypothesis. We considered 23 different values for $\Delta^2$, which are 0, 0.1, 0.2, 0.3, 0.4, 0.5, 0.6, 0.7, 0.8, 0.9, 1, 1.5, 2, 3, 5, 10, 15, 20, 25, 30, 35, 40, and 50.  Each realization was repeated 1000 times to obtain bias-squared and variance of the estimated regression parameters. Subsequently, MSEs were calculated for the least squared estimator (LSE), improved preliminary test estimator (IPT) which is based on LSE, SLE, and PSLE (READ/CHECK THESE LINES AGAIN). The responses were simulated from the following model:
\[
y_i = \sum_{i=1}^{p}x_i\beta_i + e_i
\]
where $e_i \sim N(0, \sigma^2)$ with two different values for $\sigma^2: 10^2, 20^2$. 

Secondly, the data generation setup was further modified to accommodate the number of non-zero $\beta$s in the model. In particular, we partitioned $\beta$ as $\beta=(k, q)'$ where $k$ indicates number of nonzero $\beta$s, and $q$ indicates $p-k$ zeros--a function of $\Delta^2$. To translate the above, when $p=10,$ and $k=3$, we would have $\beta=(1, 1, 1, 0, 0, 0,0,0,0,0)'$ to generate the response under null hypothesis. When $\Delta^2$ is introduced, e.g., $\Delta^2=5$, we would have $\beta=(1, 1, 1, 0,0,0,0,0,0,5)'$ to generate the response under alternative hypothesis. Clearly, inclusion of $\Delta^2$ acts as a degree of violation of the null hypothesis. As $\Delta^2$ increases, so does the degree of violation of the null hypothesis. We study the performance of the proposed estimators under varying degree of violation of null hypothesis as measured by $\Delta^2$.

Finally, the relative efficiencies were calculated using the following formula. 
\begin{equation}\label{eq:relative:efficiency}
\mbox{Relative Efficieicy (RelEff)} = \frac{\mbox{MSE}(\hat{\beta}_n^{\textrm{LE}})}{\mbox{MSE}(\hat{\beta}_n^{*})},
\end{equation}
where $\hat{\beta}^{LE}$ is the LASSO estimator, and $\hat{\beta}_n^{*}$ is one of the estimators whose relative efficiency is to be computed. As in equation (\ref{eq:relative:efficiency}), a relative efficiency greater than 1 would indicate superiority of the proposed estimator compared to the LASSO estimator. On the other hand, a relative efficiency of equal to or less than 1 would indicate that the efficiency of the estimator is at or below that of the LASSO estimator. We used R statistical software \citep{r:2014} to carry out the simulation. For obtaining LASSO estimate, {\tt glmnet()} R package \citep{glmnet:package} was used.

In the following, we discuss the results of our simulation studies. 

\subsection{Discussion of simulation results} \label{sec:simulation:discussion}
In this study, data have been generated with correlation between the $x$'s, $r=0, 0.5, 0.9$ for $n=100$, $p=10, 20, 50$, and the error variance $\sigma^2= 10^2, 20^2$. The relative efficiencies of the proposed estimators are presented in Tables~\ref{tb:rel:eff:cor0:p10:sigma10} through \ref{tb:rel:eff:cor9:p50:sigma20}.  To visually compare the results of various configurations, relative efficiency of the estimators are compared to LASSO estimator as shown in  Figures~\ref{fig:PSLE:r0} through \ref{fig:PSLE:r9}. Since positive-rule shrinkage-LASSO estimator (PSLE) outperforms all other estimators for most of the configurations, we separately compared its performance for $p=50$ and $\sigma=20$ at various correlation coefficient $r=0.0, 0.5, 0.9$, and the results are displayed in Figures \ref{fig:r0p50sigma20} through \ref{fig:r9p50sigma20}. In all of these figures, a horizontal line was drawn at 1 on the $y$-axis to facilitate the comparison among the estimators. For a given estimator, any point above this line indicates superiority of the estimator compared to the LASSO estimator in terms of relative efficiency.

The findings of simulation studies may be summarized as follows.

\begin{enumerate}[(i)]
\item The performance of LSE, LE, SLE and PSLE may be ordered as PSLE $>$ SLE $>$ LE $>$ LSE (where $>$ indicates dominance) for all $\Delta^2$. See Tables \ref{tb:rel:eff:cor0:p10:sigma10}-\ref{tb:rel:eff:cor9:p50:sigma20}. 

\item LASSO dominates over LSE uniformly.

\item At any given $p$ and $\sigma^2$, relative efficiency of the estimators is a decreasing function of both $k$ and $\Delta^2$. As $\Delta^2$ and $k$ increases, the relative efficiency decrease.

\item Gain in relative efficiency of the estimators is a decreasing function of the number of $\beta$s that are zero as indicated by $k$. Figures \ref{fig:PSLE:r0}, \ref{fig:PSLE:r5}, and \ref{fig:PSLE:r9} show that the relative efficiency of SLE and PSLE are the highest when $k=0$ (where relative efficiency is above 9). For $k=1$ the relative efficiency is around 6. The relative efficiency is around 4 for $k=$ 3 and 5. 

\item PTLE dominates IPT uniformly. However, IPT dominates the least squares estimator, which is consistent with the results found in literature. 

\item Neither PTLE nor Stein-type LASSO estimator dominate each other. 
\end{enumerate}

\section{Real data examples}\label{sec:real:data:example}

In the following, we study three real life examples. We pre-process the data sets by centering the predictor variables. We then fit linear regression models to predict the variable of interest using the available regressors. Lasso estimator (LE), improved preliminary test (IPT), Stein-type shrinkage LASSO (SLE), and Stein-type positive-rule shrinkage LASSO (PSLE) estimators are then obtained for the regression parameters. 

We evaluate the performance of the estimators by computing average cross validation error using $K=10$-fold cross validation. In cross validation, the data set is randomly divided into $K$ subsets of roughly equal size. One subset is left aside, termed as test set, while the remaining $K-1$ subsets, called training set, are used to fit the model. The fitted model is then used to predict the responses for the test data set. Finally, prediction errors are obtained by taking the squared deviation of the observed and predicted values in the test set.

In cross validation, the estimated prediction error varies across runs. Therefore, we repeat the process 1000 times, and calculate the average and standard deviation of the prediction errors. We found 1000 to be large enough number of runs to stabilize the standard deviations as no noticeable changes were observed for larger values.

\subsection{Galapagos data}


\cite{faraway:2002} analyzed the data about species diversity on the Galapagos islands. The Galapagos data contains 30 rows and seven variables. Each row represents an island, and the covariates represent various geographic measurements. The covariates are:  the number of endemic species, the area of the island (km$^2$), the highest elevation of the island (m), the distance from the nearest island (km), the distance from Santa Cruz island (km), the area of the adjacent island (km$^2$). The original data set contained missing values for some of the covariates, which have been imputed by \cite{faraway:2002} for convenience. The response variable is the number of species of tortoise found on the island. 

\begin{figure}[h!]
\centering
\includegraphics[width=4in]{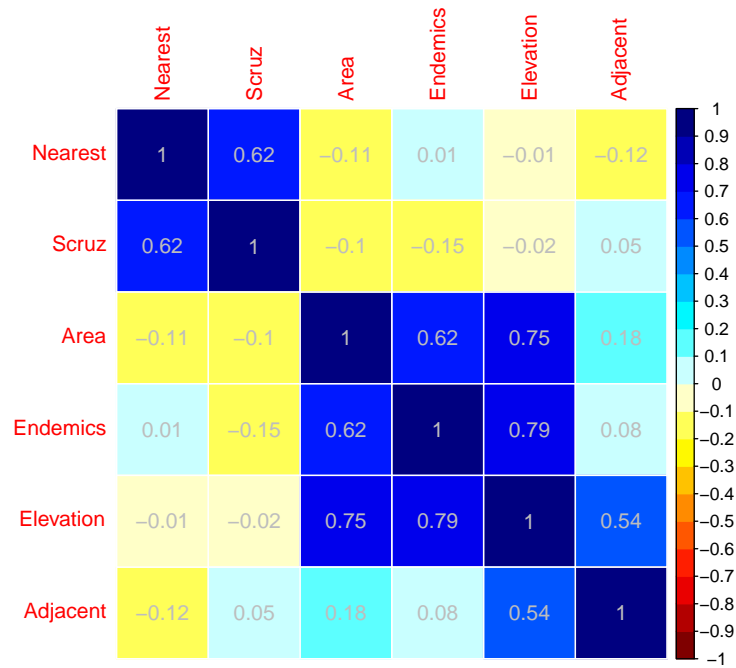} \label{fig:corr:galapagos}
\caption{Correlation matrix for Galapagos data.}
\end{figure}

The summary statistics in shown in Table~\ref{tb:summary:galapagos}. The visual correlation matrix for the centered covariates is shown in Figure~\ref{fig:corr:galapagos}.

\begin{table}[h!]
\centering
\caption{Summary statistics for the Galapagos Data. } \label{tb:summary:galapagos}
\bigskip
\begin{tabular}{ccccccc}
\hline
Min & Q1 & Median & Mean & Q3 & Max & SD (Response) \\
\hline
2.00 & 13.00 & 42.00 & 85.23 & 96.00 & 444.00 & 114.63 \\
\hline
\end{tabular}
\end{table}

\begin{figure}[h!]
\centering
\includegraphics[width=5in]{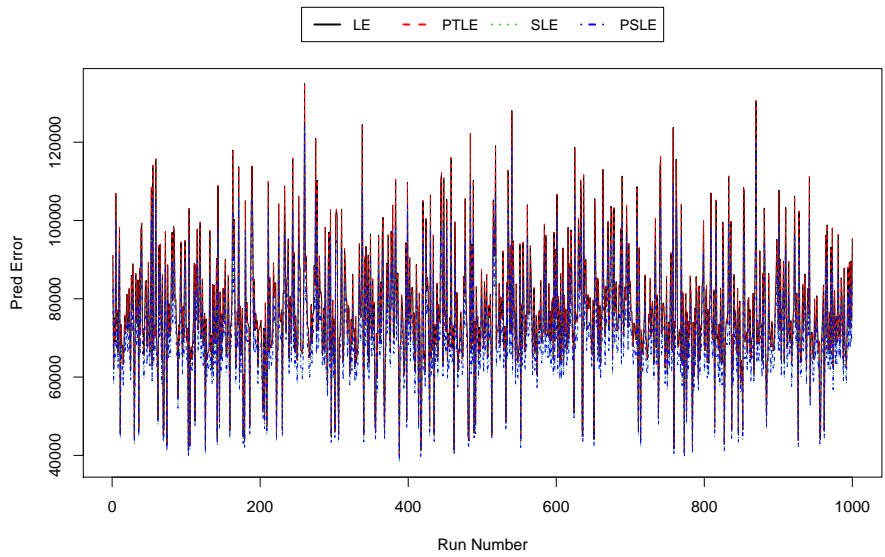} \label{fig:prederr:galapagos}
\caption{Prediction error for the proposed estimators for Galapagos data.}
\end{figure}

\begin{table}[h!]
  \centering
\caption{Average prediction errors and standard deviations for the estimators based on $K=10$-fold cross validation, repeated 1000 times for Galapagos data.}\label{tb:galapagos:est}
\medskip
\begin{tabular}{l r r} \hline
& \multicolumn{2}{c}{Bias Corrected CVE}  \\
\cline{2-3}
Estimator & Average & SD
\\
\hline
LE		& 76350.90	& 15057.76 \\
PTLE	& 76350.90	& 15057.76 \\
SLE		& 70569.66  & 13152.35 \\
PSLE	& 70569.66	& 13152.35 \\
\hline
\end{tabular}
\end{table}

Figure~\ref{fig:prederr:galapagos} displays the prediction error of the estimators for 1000 cross-validated runs. The 
average and standard deviation of the predictor errors are summarized in Table~\ref{tb:galapagos:est}. It is noted here that the prediction errors are unusually large for this data set. The reason is due to the variability present in the original data set (standard deviation for the response variable is 114.63). 

We find SLE and PSLE to be performing better than the LE and PTLE. Notably, SLE and PSLE have smaller standard deviations compared to LE and PTLE. 

\clearpage

\subsection{State data}
\cite{faraway:2002} illustrated variable selection methods using the state data set. There are 97 observations (cases) on 9 variables. The variables are: population estimate as of July 1, 1975; per capita income (1974); illiteracy (1970, percent of population); life expectancy in years (1969-71); murder and non-negligent manslaughter rate per 100,000 population (1976); percent high-school graduates (1970); mean number of days with minimum temperature 32 degrees (1931-1960) in capital or large city; and land area in square miles. We consider life expectancy as the response. 

\begin{figure}[h!]
\centering
\includegraphics[width=4in]{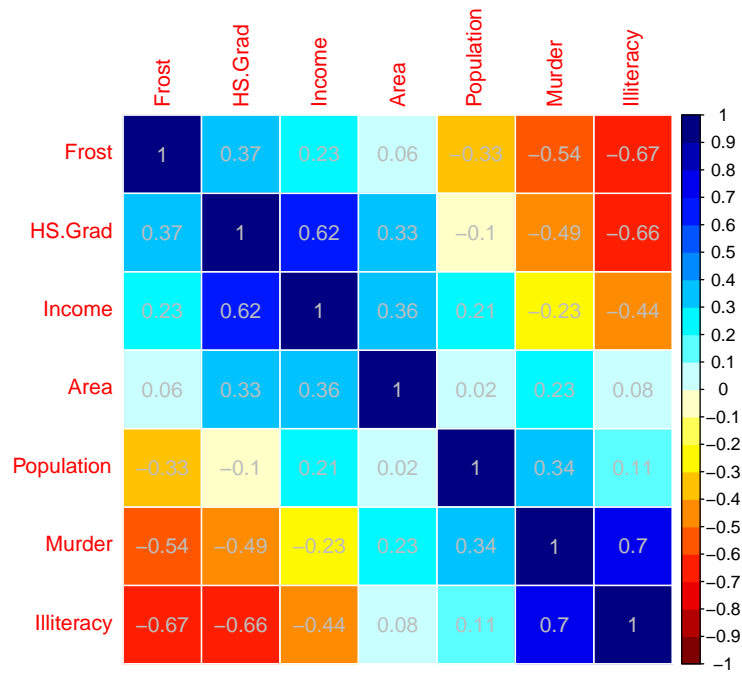} \label{fig:corr:state}
\caption{Correlation matrix for State data.}
\end{figure}

Summary statistics for this data set is given in Table \ref{tb:summary:state}. Correlation coefficients between the predictors is displayed in Figure \ref{fig:corr:state}. We notice moderate to strong correlation present between some of the predictors. In the simulation studies, we have observed that the proposed estimators perform superiorly when the correlation between the predictors is large. 

\begin{table}[h!]
\centering
\caption{Summary statistics for the State Data. } \label{tb:summary:state}
\bigskip
\begin{tabular}{ccccccc}
\hline
Min & Q1 & Median & Mean & Q3 & Max & SD (Response) \\
\hline
  67.96&   70.12&   70.68 &  70.88&   71.89&   73.60&  1.34 \\
\hline
\end{tabular}
\end{table}

\begin{table}[h!]
  \centering
\caption{Average prediction errors and standard deviations for the estimators based on $K=10$-fold cross validation, repeated 1000 times for State data.}\label{tb:state:est}
\medskip
\begin{tabular}{l rr} \hline
& \multicolumn{2}{c}{Bias Corrected CVE}  \\
\cline{2-3}
Estimator & Average & SD
\\
\hline
LE		& 5786.94 & 148.94 	\\
PTLE	& 5530.89 & 50.70 \\
SLE		& 5216.48 & 195.24 \\
PSLE	& 5552.93 & 81.89 \\
\hline
\end{tabular}
\end{table}

\begin{figure}[h!]
\centering
\includegraphics[width=5in]{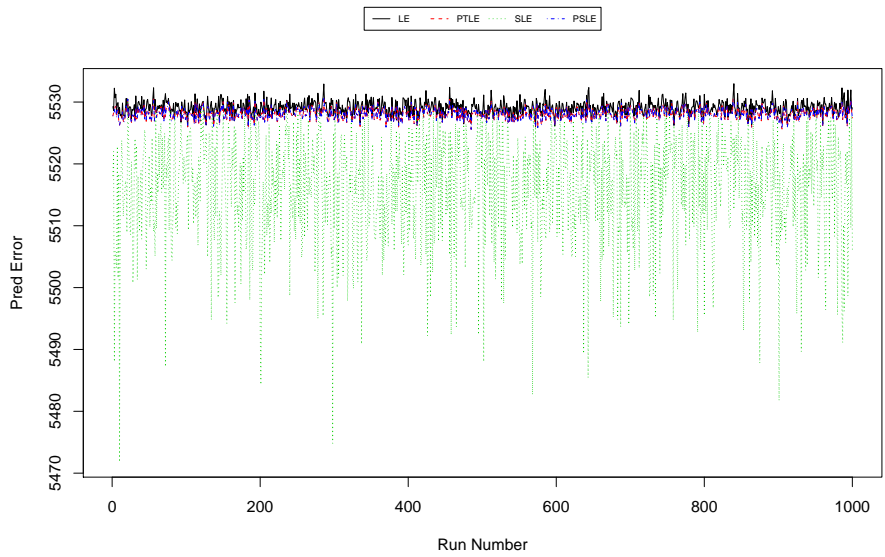} \label{fig:prederr:state}
\caption{Prediction error of the proposed estimators for State data.}
\end{figure}

Table  \ref{tb:state:est} gives averages and standard deviations of the predictor errors for the estimators. For this data, SLE has the smallest average prediction error followed by PSLE. LE has the largest average prediction error. Figure \ref{fig:prederr:state} shows the average prediction errors for the estimators visually, which demonstrates smaller yet highly variable prediction errors for the SLE. 

\clearpage

\subsection{Longley data}

\begin{figure}[h!]
\centering
\includegraphics[width=4.5in]{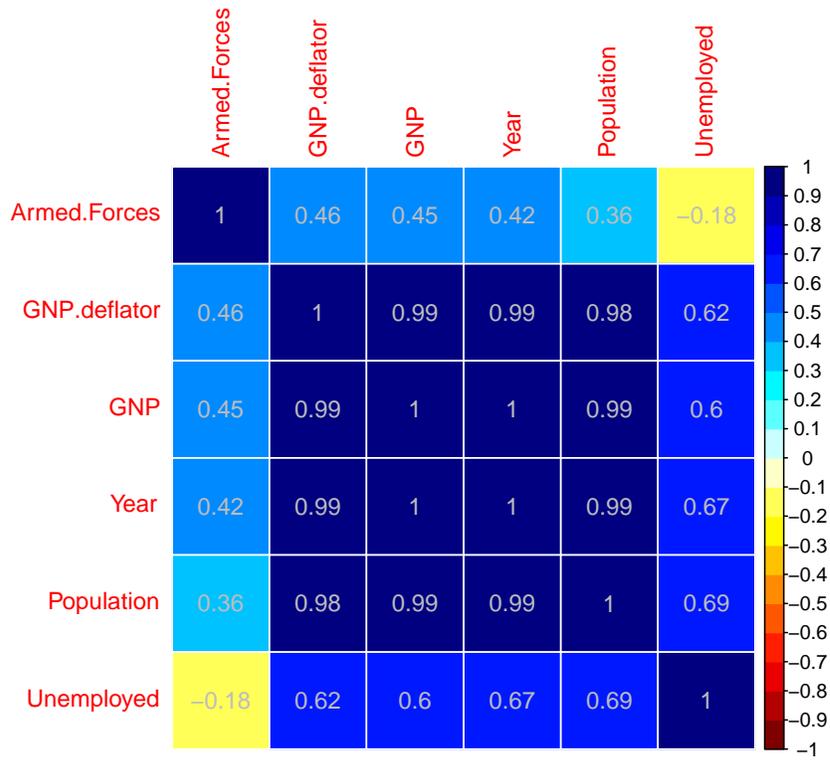} \label{fig:corr:long}
\caption{Correlation matrix for Longley data.}
\end{figure}

Longley data set \citep{longley:1967} is a popular macroeconomic data set and is widely used to demonstrate the application of collinear regression. The data has 16 observations on seven variables. The predictors are GNP implicit price deflator, gross national product, the number of people unemployed, number of people in the armed forces, noninstitutionalized population aged 14 and older, the year (time). The response variable is the number of people employed. Table \ref{tb:summary:long} shows data summary while visual correlation matrix for the predictors is shown in Figure \ref{fig:corr:long}. Note that the predictors are highly correlated. Table \ref{tb:long:est} shows SLE to be the best estimator in terms of prediction errors followed by PTLE and PSLE. LE has the largest average prediction error. Although SLE has the smallest prediction error, it has the largest variability. Of the four estimators, PTLE has the smallest variability. The plot of 1000 prediction errors in Figure~\ref{fig:corr:long} demonstrates SLE as the best performing estimator.

\begin{table}[h!]
\centering
\caption{Summary statistics for the longley Data. } \label{tb:summary:long}
\bigskip
\begin{tabular}{ccccccc}
\hline
Min & Q1 & Median & Mean & Q3 & Max & SD (Response) \\
\hline
60.17 &  62.71 &  65.50 &  65.32 &  68.29 &  70.55 & 3.51 \\
\hline
\end{tabular}
\end{table}

\begin{table}[h!]
  \centering
\caption{Average prediction errors and standard deviations for the estimators based on $K=10$-fold cross validation, repeated 1000 times for Longley data.}\label{tb:long:est}
\medskip
\begin{tabular}{l rr} \hline
& \multicolumn{2}{c}{Bias Corrected CVE}  \\
\cline{2-3}
Estimator & Average & SD
\\
\hline
LE		& 5182.43 & 99.16 \\
PTLE	& 4751.56 & 5.97 \\
SLE		& 3163.05 & 337.86 \\
PSLE	& 4759.06 & 24.25 \\
\hline
\end{tabular}
\end{table}

\begin{figure}[h!]
\centering
\includegraphics[width=5in]{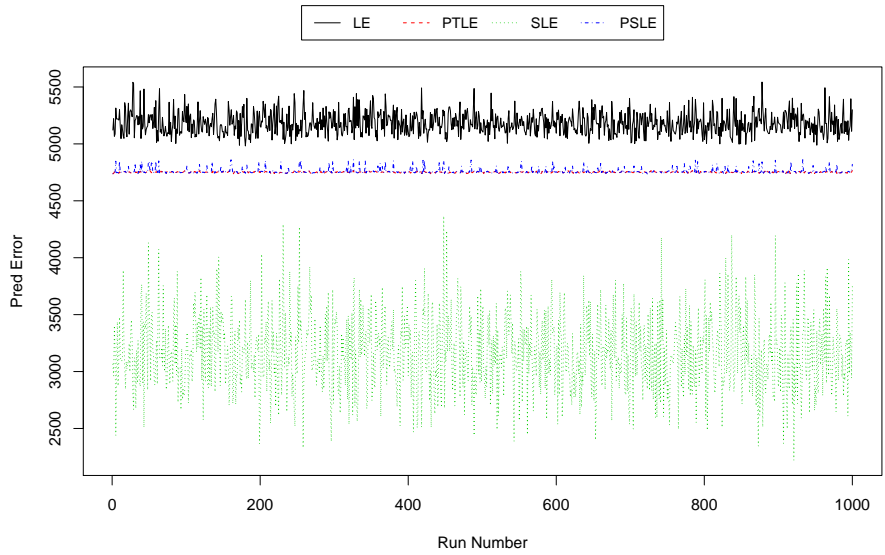} \label{fig:prederr:long}
\caption{Prediction error of the proposed estimators for Longley data.}
\end{figure}

\clearpage

\section{Discussion and conclusion}
\begin{figure}[h!]
\includegraphics[width=6in]{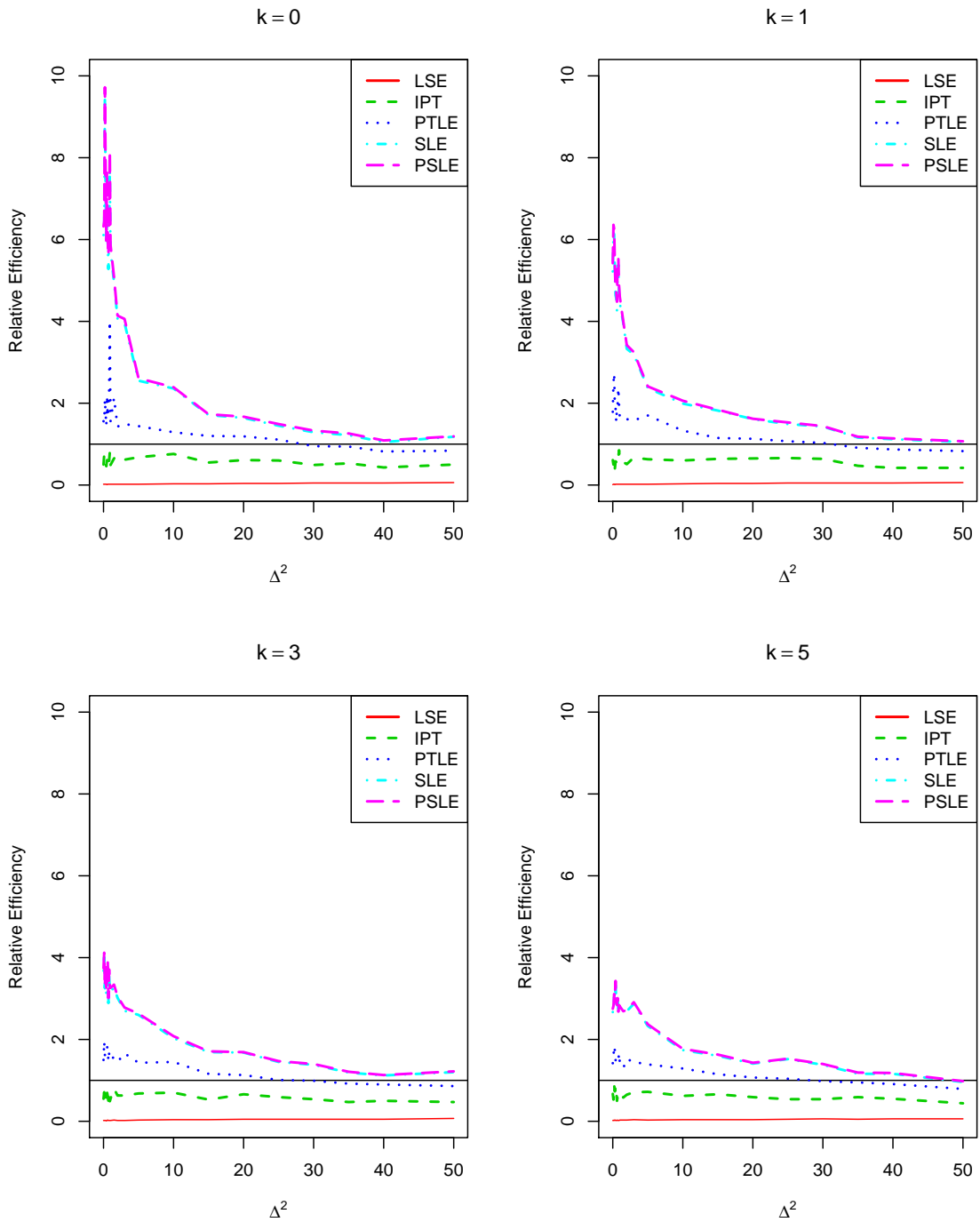} \label{fig:r0p50sigma20}
\caption{Relative efficiencies of LSE, IPT, PTLE, SLE, PSLE, compared to LASSO estimator (LE) when  $r=0$, $p=50$, and $\sigma=20$.}
\end{figure}

\begin{figure}[h!]
\includegraphics[width=6in]{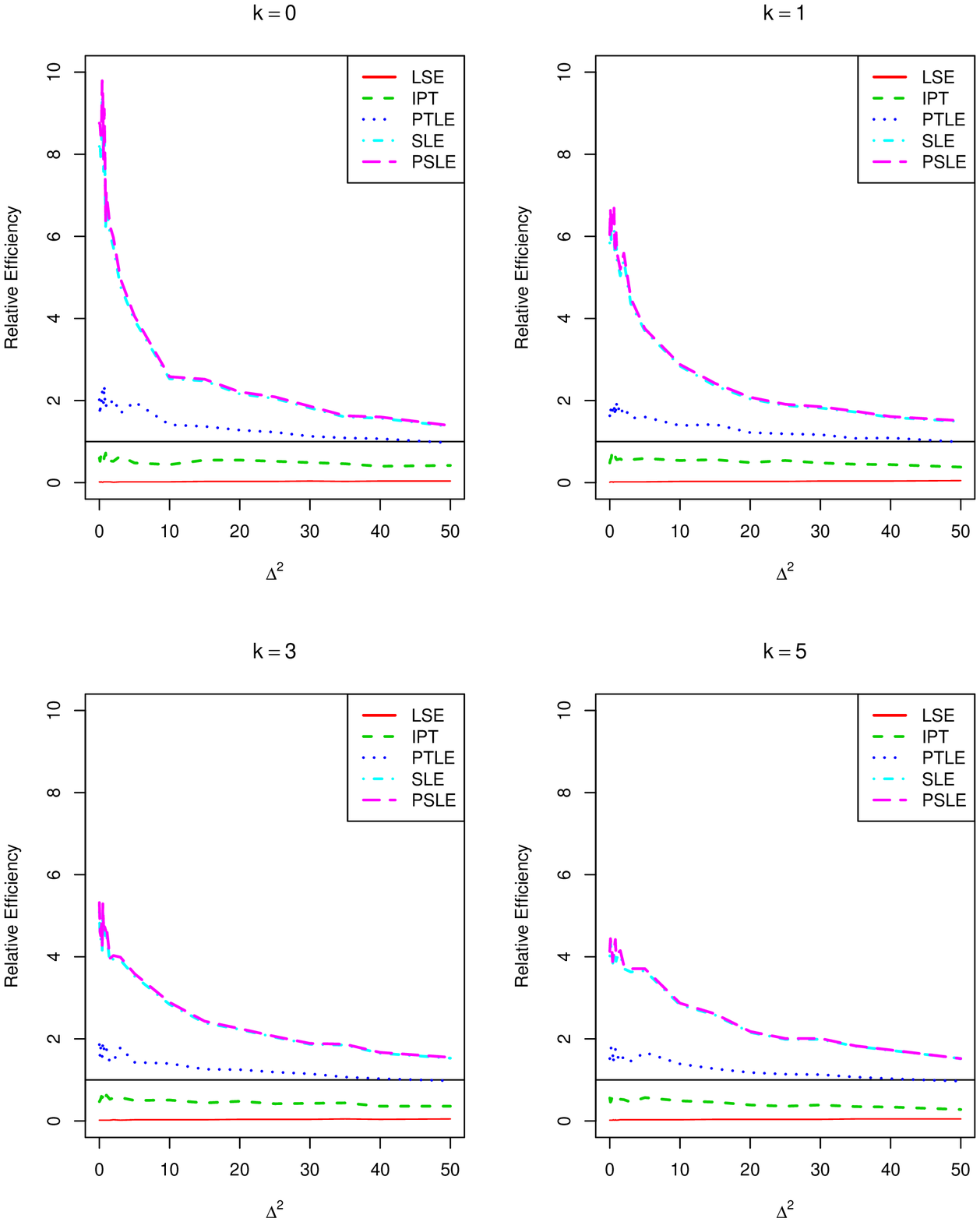} \label{fig:r5p50sigma20}
\caption{Relative efficiencies of LSE, IPT, PTLE, SLE, PSLE, compared to LASSO estimator (LE) when  $r=0.5$, $p=50$, and $\sigma=20$.}
\end{figure}

\begin{figure}[h!]
\includegraphics[width=6in]{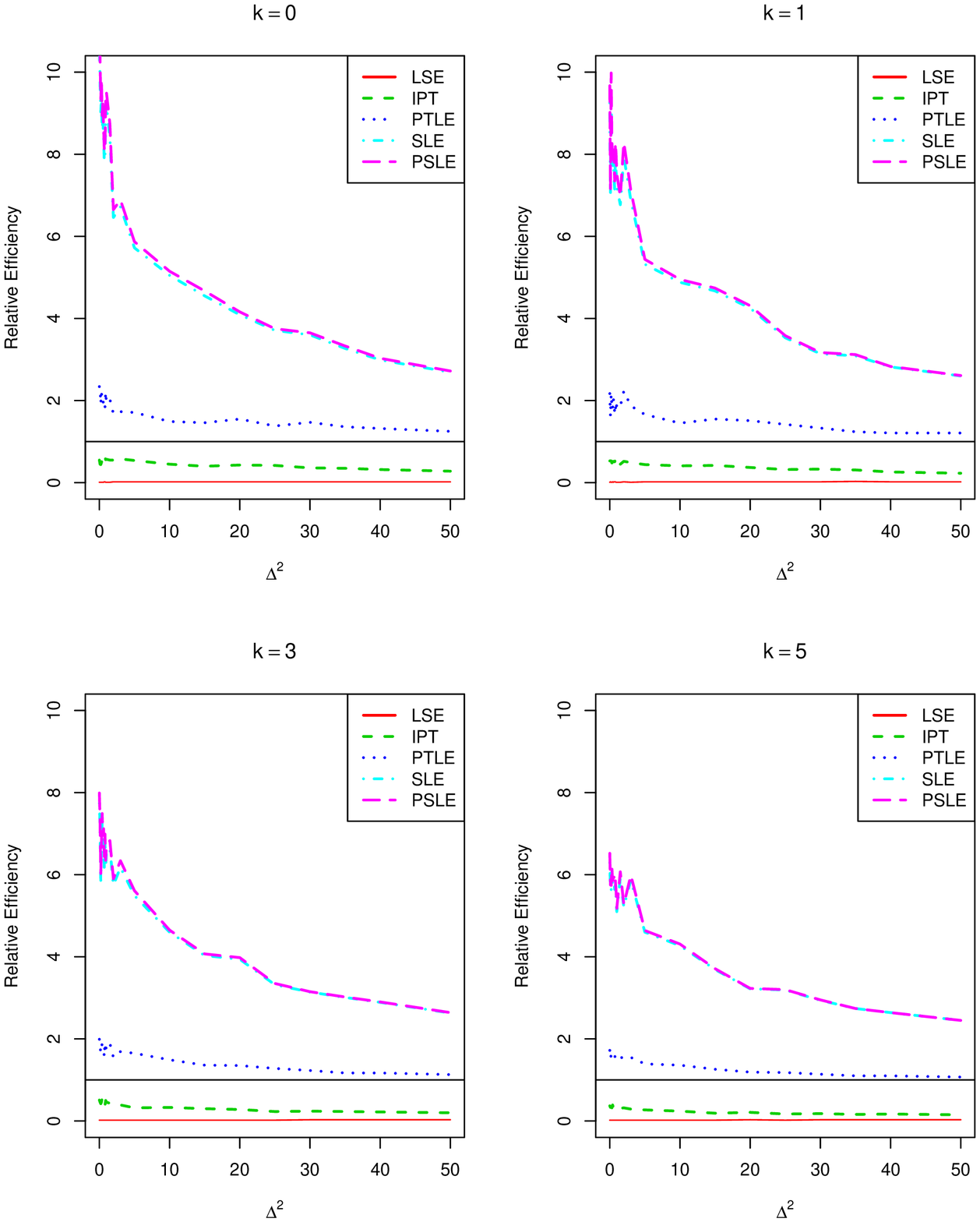} \label{fig:r9p50sigma20}
\caption{Relative efficiencies of LSE, IPT, PTLE, SLE, PSLE, compared to LASSO estimator (LE) when  $r=0.9$, $p=50$, and $\sigma=20$.}
\end{figure}

\begin{figure}[h!]
\includegraphics[width=6in]{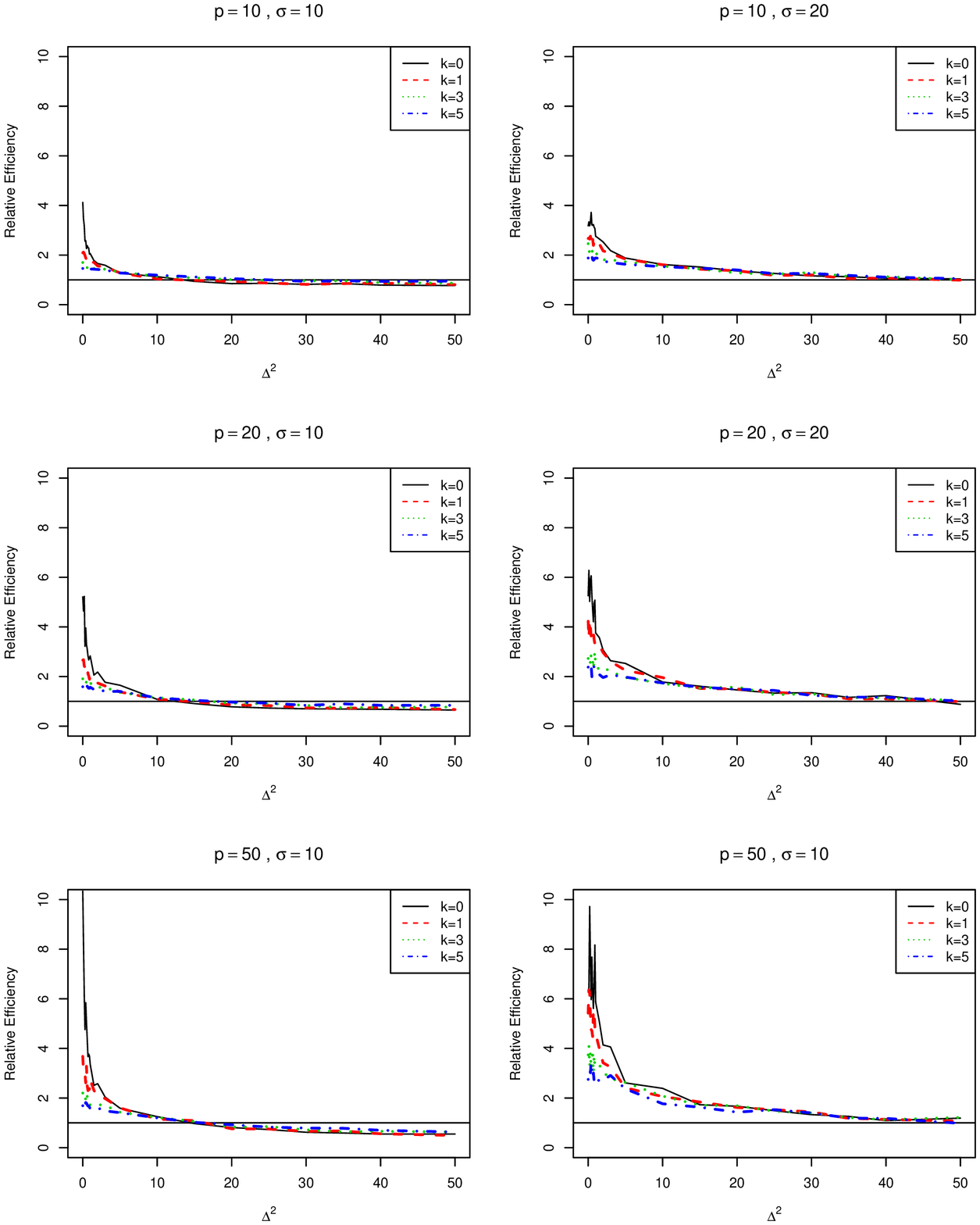}\label{fig:PSLE:r0}
\caption{Relative efficiency of positive rule shrinkage LASSO estimator (PSLE) compared to LASSO estimator (LE) for various $k$ when  $r=0$.}
\end{figure}

\begin{figure}[h!]
\includegraphics[width=6in]{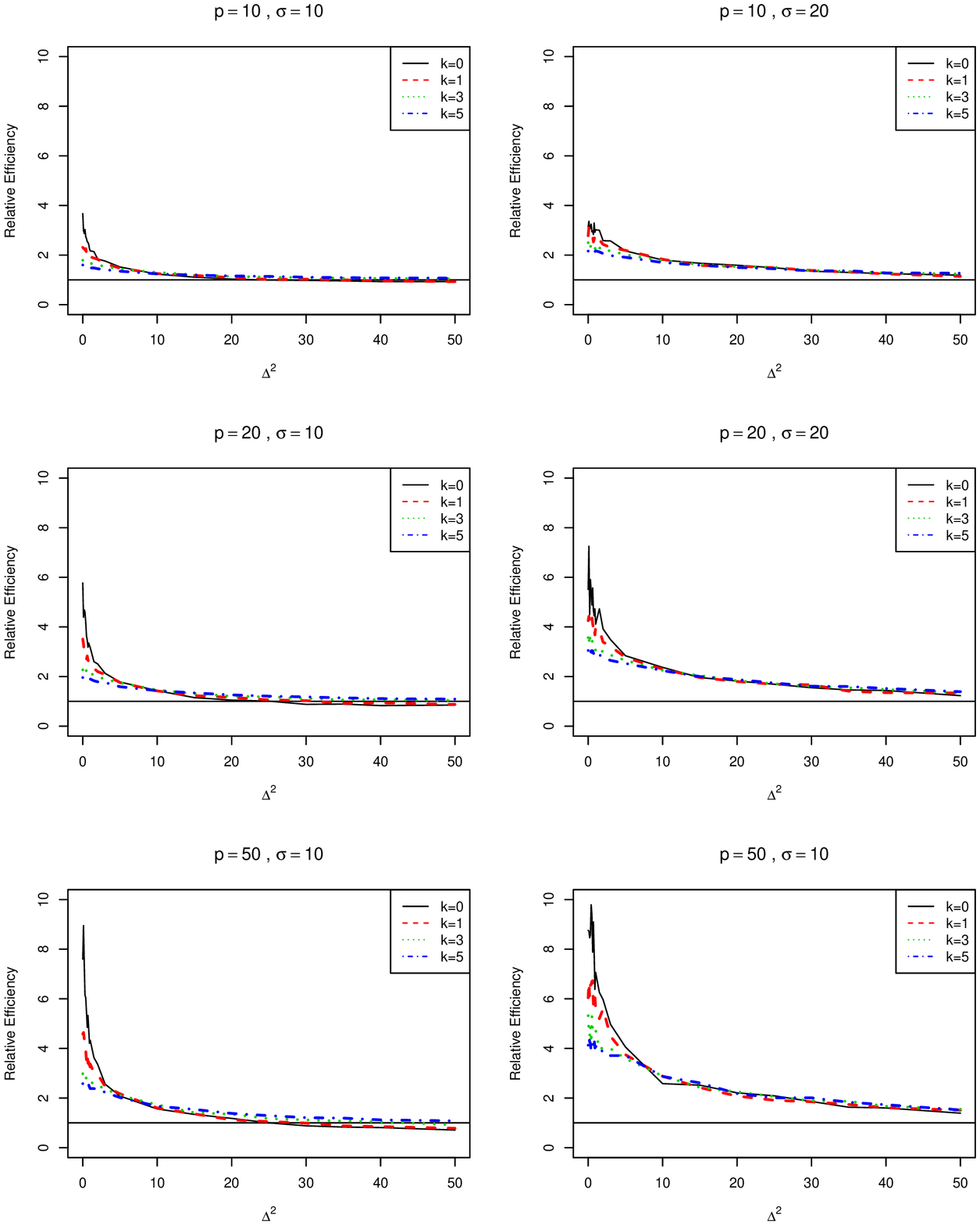} \label{fig:PSLE:r5}
\caption{Relative efficiency of positive rule shrinkage LASSO estimator (PSLE) compared to LASSO estimator (LE) for various $k$ when  $r=0.5$.}
\end{figure}

\begin{figure}[h!]
\includegraphics[width=6in]{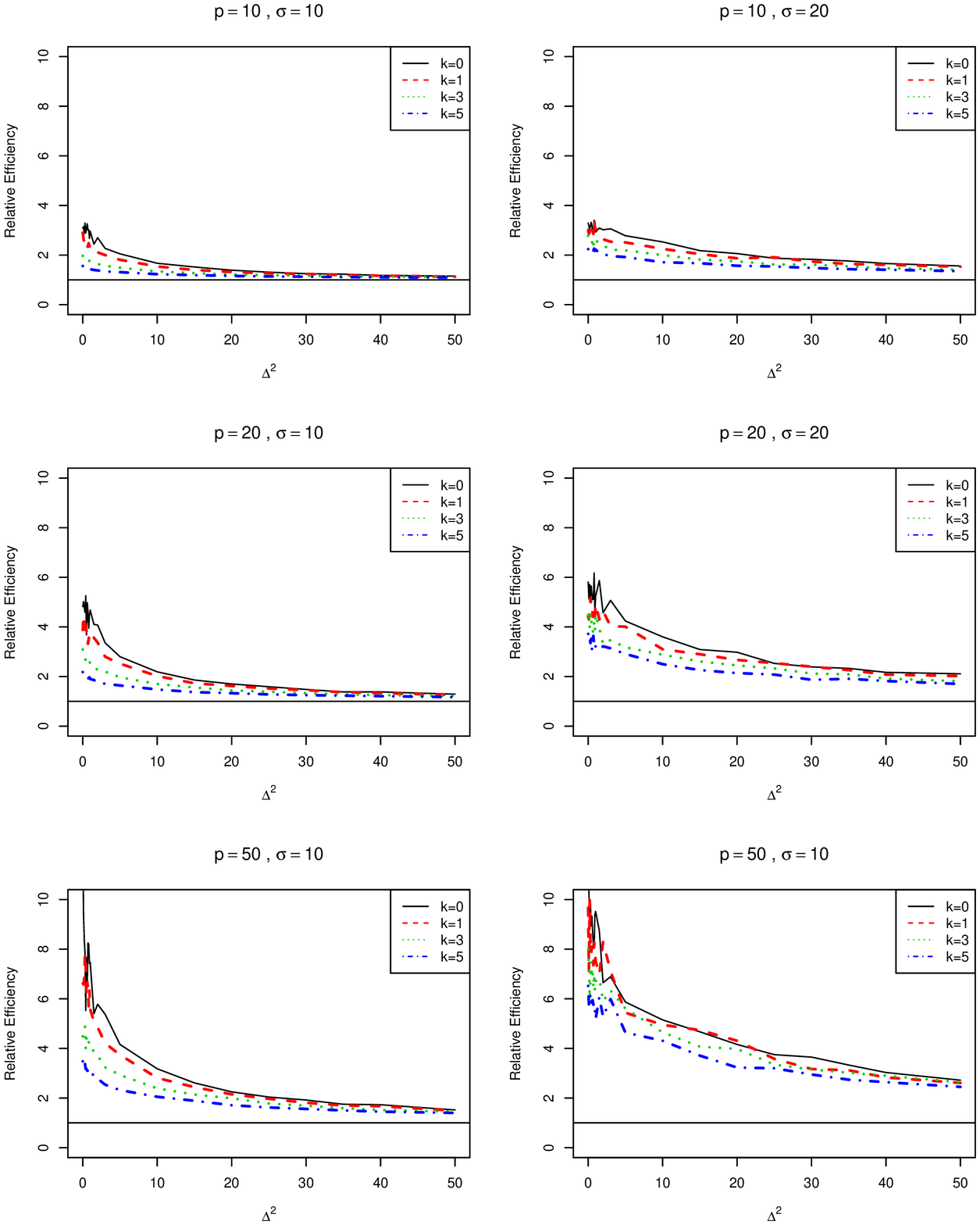} \label{fig:PSLE:r9}
\caption{Relative efficiency of positive rule shrinkage LASSO estimator (PSLE) compared to LASSO estimator (LE) for various $k$ when  $r=0.9$.}
\end{figure}

In this paper, we proposed improved estimation technique for the LASSO estimator based on Stein-rule. In particular, we used LASSO estimator to obtain preliminary test LASSO estimator (PTLE), Stein-type shrinkage LASSO estimator (SLE), and Stein-type positive-rule shrinkage LASSO estimator (PSLE). We studied the performance of the proposed estimators under full-model hypothesis when the parameter space is small relative to sample size. Simulation studies have been performed to compare the estimators for various configurations of parameter sizes ($p$), correlation coefficient between the predictors ($r$), and the error variance ($\sigma^2$). We used relative-MSE criterion to compare the proposed estimators with the classical LASSO estimator. We varied the number of non-zero $\beta$s, and evaluated the performance of the estimators under varying degree of model misspecification as guided by $\Delta^2$. We have provided relative efficiencies of the proposed estimators compared to classical LASSO estimator in Tables~\ref{tb:rel:eff:cor0:p10:sigma10} through \ref{tb:rel:eff:cor9:p50:sigma20} as well as graphically for selected configurations as displayed in Figures \ref{fig:r0p50sigma20} through \ref{fig:PSLE:r9}. 

The simulation results demonstrate that the classical LASSO dominates the LSE uniformly while PSLE has the smallest MSE among the proposed estimators. In particular, PSLE uniformly dominates classical LASSO estimator when the error variance is large ($\sigma^2=20^2$, in our setup). Also, neither PTLE nor SLE dominates one another. Relative efficiency of the proposed estimators increases when there are more near-zero parameters present in the model. Performance of the estimators decrease as we deviate from the null model. These results are consistent with the properties of traditional preliminary test and Stein-type estimators found in the literature. 

We have presented three real life data examples to demonstrate the use of proposed estimators. Average and standard deviations of the prediction errors based on the LE, PTLE, SLE, and PSLE have been obtained and compared. The proposed estimators outperform LASSO estimator in both average prediction error and standard deviation criteria. While PSLE dominated the other proposed estimators in the simulation experiments, the dominance picture was not obvious for the real data examples. As such we conclude that neither of the PTLE, SLE or PSLE may outperform each other in all real life applications. We tried some other data sets where the correlation between the predictors are low to moderately strong, and found that LASSO as well as the proposed estimators perform equally in those cases. Therefore, we conclude that the improved LASSO estimators would find their applications for data sets with moderate to strong correlations among the predictors.

\clearpage

\clearpage

\clearpage

\pagebreak

\pagebreak

\pagebreak

\begin{landscape}
\begin{table}[!h]
\centering
\caption{Relative efficiencies of IPT, PTLE, SLE, and PSLE with respect to LASSO estimator when $n=100$, $r=0,$ $p=10$, $\sigma=10$ for different values of $k$.} \label{tb:rel:eff:cor0:p10:sigma10}
\bigskip
\begin{footnotesize}
\begin{tabular}{r|rrrrr|rrrrr|rrrrr|rrrrr}
\hline
 & \multicolumn{5}{c|}{$k=0$} & \multicolumn{5}{c|}{$k=1$} & \multicolumn{5}{c|}{$k=3$}
 & \multicolumn{5}{c}{$k=5$} \\
 $\Delta^2$ & LSE & IPT & PTLE & SLE & PSLE & LSE & IPT & PTLE & SLE & PSLE & 
 LSE & IPT & PTLE & SLE & PSLE& LSE & IPT & PTLE & SLE & PSLE \\
 \hline
 0.0&0.10&1.84&2.50&4.06&4.12&0.17&1.24&1.44&2.08&2.09&0.28&1.24&1.21&1.70 & 1.70 &0.38 &1.02 &1.03 &1.45 &1.46\\
 0.1&0.10&2.82&3.04&3.41&3.46&0.16&1.36&1.68&2.10&2.12&0.26&1.06&1.16&1.53 & 1.53 &0.37 &1.00 &1.04 &1.41 &1.42\\
 0.2&0.12&1.79&2.37&3.17&3.18&0.17&1.33&1.59&2.04&2.06&0.27&1.07&1.13&1.55 & 1.55 &0.40 &1.06 &1.04 &1.43 &1.43\\
 0.3&0.11&1.67&1.62&2.56&2.56&0.18&1.35&1.36&1.97&1.98&0.31&1.10&1.09&1.57 & 1.58 &0.39 &1.02 &1.03 &1.41 &1.42\\
 0.4&0.11&1.52&1.78&2.57&2.58&0.19&1.32&1.28&1.88&1.89&0.28&1.07&1.12&1.50 & 1.50 &0.40 &1.06 &1.04 &1.45 &1.46\\
 0.5&0.14&1.50&1.53&2.24&2.26&0.18&1.21&1.34&1.78&1.80&0.31&1.10&1.12&1.52 & 1.53 &0.43 &1.06 &1.03 &1.44 &1.44\\
\hline
\end{tabular}
\end{footnotesize}
\end{table}
\end{landscape}

\begin{landscape}
\begin{table}[!h]
\centering
\caption{Relative efficiencies of IPT, PTLE, SLE, and PSLE with respect to LASSO estimator when $n=100$, $r=0,$ $p=10$, $\sigma=20$ for different values of $k$.} \label{tb:rel:eff:cor0:p10:sigma20}
\bigskip
\begin{footnotesize}
\begin{tabular}{r|rrrrr|rrrrr|rrrrr|rrrrr}
\hline
 & \multicolumn{5}{c|}{$k=0$} & \multicolumn{5}{c|}{$k=1$} & \multicolumn{5}{c|}{$k=3$}
 & \multicolumn{5}{c}{$k=5$} \\
 $\Delta^2$ & LSE & IPT & PTLE & SLE & PSLE & LSE & IPT & PTLE & SLE & PSLE & 
 LSE & IPT & PTLE & SLE & PSLE& LSE & IPT & PTLE & SLE & PSLE \\
 \hline
 0.0&0.09& 1.43& 1.89& 3.17& 3.18&0.11& 1.36& 1.69& 2.66& 2.67& 0.14& 1.39& 1.90& 2.45& 2.46& 0.17& 1.43& 1.40& 1.88& 1.88\\
 0.1&0.11& 1.68& 1.80& 3.31& 3.34&0.13& 1.57& 1.60& 2.68& 2.69& 0.14& 1.32& 1.50& 2.09& 2.09& 0.17& 1.38& 1.45& 2.00& 2.00\\
 0.2&0.12& 1.90& 1.95& 3.17& 3.18&0.11& 1.70& 1.85& 2.62& 2.64& 0.13& 1.41& 1.64& 2.18& 2.19& 0.19& 1.35& 1.32& 1.90& 1.90\\
 0.3&0.11& 1.81& 2.17& 3.30& 3.33&0.12& 1.72& 2.28& 2.75& 2.76& 0.15& 1.26& 1.40& 2.09& 2.10& 0.18& 1.40& 1.39& 1.93& 1.93\\
 0.4&0.09& 2.39& 2.50& 3.69& 3.72&0.14& 1.31& 1.44& 2.52& 2.52& 0.17& 1.35& 1.66& 2.18& 2.19& 0.18& 1.21& 1.31& 1.90& 1.91\\
 0.5&0.13& 2.16& 2.23& 3.19& 3.21&0.13& 1.64& 1.84& 2.75& 2.77& 0.17& 1.50& 1.52& 2.26& 2.26& 0.17& 1.23& 1.43& 1.96& 1.96\\
\hline
\end{tabular}
\end{footnotesize}
\end{table}
\end{landscape}

\begin{landscape}
\begin{table}[!h]
\centering
\caption{Relative efficiencies of IPT, PTLE, SLE, and PSLE with respect to LASSO estimator when $n=100$, $r=0,$ $p=20$, $\sigma=10$ for different values of $k$.} \label{tb:rel:eff:cor0:p20:sigma10}
\bigskip
\begin{footnotesize}
\begin{tabular}{r|rrrrr|rrrrr|rrrrr|rrrrr}
\hline
 & \multicolumn{5}{c|}{$k=0$} & \multicolumn{5}{c|}{$k=1$} & \multicolumn{5}{c|}{$k=3$}
 & \multicolumn{5}{c}{$k=5$} \\
 $\Delta^2$ & LSE & IPT & PTLE & SLE & PSLE & LSE & IPT & PTLE & SLE & PSLE & 
 LSE & IPT & PTLE & SLE & PSLE& LSE & IPT & PTLE & SLE & PSLE \\
 \hline
 0.0& 0.06& 1.33& 1.81& 5.18& 5.21& 0.07& 1.41& 1.73& 2.63& 2.66& 0.14& 1.05& 1.22& 1.90& 1.91& 0.17& 1.07& 1.12& 1.58& 1.59\\
 0.1& 0.05& 1.40& 1.96& 4.58& 4.64& 0.08& 1.37& 1.70& 2.65& 2.68& 0.13& 0.95& 1.18& 1.80& 1.80& 0.16& 0.96& 1.10& 1.49& 1.49\\
 0.2& 0.06& 2.25& 2.94& 5.17& 5.23& 0.08& 1.12& 1.37& 2.41& 2.43& 0.13& 1.04& 1.21& 1.71& 1.72& 0.18& 0.96& 1.07& 1.66& 1.66\\
 0.3& 0.07& 1.16& 1.57& 3.17& 3.22& 0.08& 1.23& 1.37& 2.35& 2.36& 0.15& 1.09& 1.14& 1.84& 1.85& 0.19& 1.04& 1.14& 1.61& 1.61\\
 0.4& 0.07& 1.87& 2.35& 3.93& 3.96& 0.10& 1.22& 1.42& 2.47& 2.49& 0.13& 1.09& 1.24& 1.79& 1.80& 0.18& 0.95& 1.07& 1.59& 1.60\\
 0.5& 0.07& 1.46& 1.56& 3.25& 3.29& 0.10& 1.10& 1.36& 2.39& 2.40& 0.15& 1.08& 1.18& 1.75& 1.76& 0.19& 1.00& 1.08& 1.60& 1.61\\
\hline
\end{tabular}
\end{footnotesize}
\end{table}
\end{landscape}

\begin{landscape}
\begin{table}[!h]
\centering
\caption{Relative efficiencies of IPT, PTLE, SLE, and PSLE with respect to LASSO estimator when $n=100$, $r=0,$ $p=20$, $\sigma=20$ for different values of $k$.} \label{tb:rel:eff:cor0:p20:sigma20}
\bigskip
\begin{footnotesize}
\begin{tabular}{r|rrrrr|rrrrr|rrrrr|rrrrr}
\hline
 & \multicolumn{5}{c|}{$k=0$} & \multicolumn{5}{c|}{$k=1$} & \multicolumn{5}{c|}{$k=3$}
 & \multicolumn{5}{c}{$k=5$} \\
 $\Delta^2$ & LSE & IPT & PTLE & SLE & PSLE & LSE & IPT & PTLE & SLE & PSLE & 
 LSE & IPT & PTLE & SLE & PSLE& LSE & IPT & PTLE & SLE & PSLE \\
 \hline
 0.0&0.07& 1.75& 2.13& 5.17& 5.26& 0.06& 1.09& 1.86& 4.01& 4.23& 0.06& 1.40& 1.56& 2.69& 2.73& 0.08& 1.05& 1.42& 2.36& 2.38\\
 0.1&0.04& 1.48& 2.53& 5.89& 6.28& 0.06& 1.59& 2.09& 3.55& 3.56& 0.06& 1.65& 1.92& 2.71& 2.72& 0.09& 1.26& 1.30& 2.39& 2.42\\
 0.2&0.06& 1.57& 2.04& 4.99& 5.03& 0.05& 1.39& 2.26& 4.20& 4.27& 0.07& 1.72& 1.83& 2.80& 2.84& 0.08& 1.06& 1.35& 2.35& 2.37\\
 0.3&0.04& 1.39& 2.53& 5.85& 5.89& 0.07& 1.42& 1.64& 3.92& 3.98& 0.06& 1.24& 1.71& 2.48& 2.49& 0.09& 1.37& 1.60& 2.48& 2.49\\
 0.4&0.05& 1.51& 2.63& 5.83& 6.06& 0.07& 1.52& 1.87& 3.60& 3.61& 0.07& 1.17& 1.67& 2.88& 2.91& 0.08& 1.12& 1.57& 2.45& 2.46\\
 0.5&0.06& 1.75& 2.42& 5.23& 5.28& 0.06& 1.66& 1.83& 3.65& 3.72& 0.07& 1.46& 1.77& 2.84& 2.86& 0.07& 1.06& 1.35& 1.98& 1.99\\
\hline
\end{tabular}
\end{footnotesize}
\end{table}
\end{landscape}

\begin{landscape}
\begin{table}[!h]
\centering
\caption{Relative efficiencies of IPT, PTLE, SLE, and PSLE with respect to LASSO estimator when $n=100$, $r=0,$ $p=50$, $\sigma=10$ for different values of $k$.} \label{tb:rel:eff:cor0:p50:sigma10}
\bigskip
\begin{footnotesize}
\begin{tabular}{r|rrrrr|rrrrr|rrrrr|rrrrr}
\hline
 & \multicolumn{5}{c|}{$k=0$} & \multicolumn{5}{c|}{$k=1$} & \multicolumn{5}{c|}{$k=3$}
 & \multicolumn{5}{c}{$k=5$} \\
 $\Delta^2$ & LSE & IPT & PTLE & SLE & PSLE & LSE & IPT & PTLE & SLE & PSLE & 
 LSE & IPT & PTLE & SLE & PSLE& LSE & IPT & PTLE & SLE & PSLE \\
 \hline
 0.0& 0.02& 0.57& 2.24& 9.33& 10.34& 0.02& 0.71& 1.52& 3.61& 3.68& 0.04& 0.70& 1.39& 2.18& 2.20& 0.04& 0.66& 1.20& 1.68& 1.69\\
 0.1& 0.02& 0.62& 1.94& 7.52&  8.26& 0.02& 0.58& 1.44& 3.04& 3.07& 0.04& 0.68& 1.27& 2.13& 2.15& 0.05& 0.71& 1.23& 1.86& 1.87\\
 0.2& 0.02& 0.58& 1.63& 5.86&  6.17& 0.02& 0.65& 1.47& 2.87& 2.99& 0.03& 0.71& 1.32& 1.94& 1.97& 0.05& 0.75& 1.18& 1.73& 1.75\\
 0.3& 0.02& 0.57& 1.75& 4.63&  4.76& 0.03& 0.57& 1.44& 2.93& 2.97& 0.03& 0.66& 1.29& 1.86& 1.88& 0.05& 0.64& 1.16& 1.71& 1.72\\
 0.4& 0.02& 0.94& 2.03& 5.59&  5.84& 0.02& 0.49& 1.46& 2.57& 2.61& 0.04& 0.59& 1.34& 1.98& 2.00& 0.05& 0.76& 1.18& 1.82& 1.82\\
 0.5& 0.02& 0.78& 2.19& 4.93&  5.05& 0.03& 0.92& 1.49& 3.23& 3.28& 0.03& 0.66& 1.24& 1.82& 1.84& 0.05& 0.76& 1.20& 1.77& 1.79\\
\hline
\end{tabular}
\end{footnotesize}
\end{table}
\end{landscape}

\begin{landscape}
\begin{table}[!h]
\centering
\caption{Relative efficiencies of IPT, PTLE, SLE, and PSLE with respect to LASSO estimator when $n=100$, $r=0,$ $p=50$, $\sigma=20$ for different values of $k$.} \label{tb:rel:eff:cor0:p50:sigma20}
\bigskip
\begin{footnotesize}
\begin{tabular}{r|rrrrr|rrrrr|rrrrr|rrrrr}
\hline
 & \multicolumn{5}{c|}{$k=0$} & \multicolumn{5}{c|}{$k=1$} & \multicolumn{5}{c|}{$k=3$}
 & \multicolumn{5}{c}{$k=5$} \\
 $\Delta^2$ & LSE & IPT & PTLE & SLE & PSLE & LSE & IPT & PTLE & SLE & PSLE &  
 LSE & IPT & PTLE & SLE & PSLE& LSE & IPT & PTLE & SLE & PSLE \\
 \hline
 0.0& 0.02& 0.50& 1.56& 6.11& 6.32& 0.01& 0.61& 1.79& 5.22& 5.42& 0.02& 0.55& 1.50& 3.68& 3.74& 0.02& 0.67& 1.42& 2.67& 2.75\\
 0.1& 0.02& 0.74& 1.53& 6.18& 6.45& 0.02& 0.50& 2.50& 6.17& 6.36& 0.02& 0.71& 1.95& 4.01& 4.12& 0.02& 0.53& 1.59& 2.75& 2.82\\
 0.2& 0.02& 0.74& 2.09& 9.49& 9.72& 0.02& 0.65& 2.73& 6.12& 6.26& 0.02& 0.47& 1.61& 3.25& 3.29& 0.03& 0.88& 1.84& 3.04& 3.11\\
 0.3& 0.02& 0.67& 1.76& 7.36& 7.74& 0.01& 0.40& 1.90& 5.15& 5.38& 0.02& 0.68& 1.69& 3.72& 3.84& 0.02& 0.73& 1.73& 2.94& 3.02\\
 0.4& 0.02& 0.46& 1.45& 5.86& 5.97& 0.02& 0.58& 1.88& 4.64& 4.76& 0.01& 0.49& 1.73& 3.10& 3.18& 0.03& 0.87& 1.49& 3.38& 3.43\\
 0.5& 0.01& 0.61& 2.19& 6.87& 7.68& 0.02& 0.56& 1.54& 4.56& 4.71& 0.02& 0.71& 1.83& 3.33& 3.44& 0.02& 0.74& 1.56& 2.84& 2.90\\
\hline
\end{tabular}
\end{footnotesize}
\end{table}
\end{landscape}

\begin{landscape}
\begin{table}[!h]
\centering
\caption{Relative efficiencies of IPT, PTLE, SLE, and PSLE with respect to LASSO estimator when $n=100$, $r=0.5,$ $p=10$, $\sigma=10$ for different values of $k$.} \label{tb:rel:eff:cor5:p10:sigma10}
\bigskip
\begin{footnotesize}
\begin{tabular}{r|rrrrr|rrrrr|rrrrr|rrrrr}
\hline
 & \multicolumn{5}{c|}{$k=0$} & \multicolumn{5}{c|}{$k=1$} & \multicolumn{5}{c|}{$k=3$}
 & \multicolumn{5}{c}{$k=5$} \\
 $\Delta^2$ & LSE & IPT & PTLE & SLE & PSLE & LSE & IPT & PTLE & SLE & PSLE & 
 LSE & IPT & PTLE & SLE & PSLE& LSE & IPT & PTLE & SLE & PSLE \\
 \hline
 0.0& 0.08& 1.98& 2.23& 3.66& 3.67& 0.14& 1.27& 1.54& 2.30& 2.31& 0.22& 0.83& 1.12& 1.79& 1.79& 0.32& 0.67& 1.01& 1.60& 1.60\\
 0.1& 0.09& 1.49& 1.86& 2.94& 2.99& 0.14& 1.20& 1.43& 2.26& 2.27& 0.24& 0.83& 1.11& 1.81& 1.81& 0.33& 0.64& 1.01& 1.55& 1.55\\
 0.2& 0.10& 1.75& 1.90& 2.86& 2.87& 0.13& 1.16& 1.47& 2.25& 2.26& 0.24& 0.79& 1.10& 1.76& 1.77& 0.31& 0.63& 1.00& 1.56& 1.56\\
 0.3& 0.10& 1.73& 2.24& 3.02& 3.03& 0.15& 1.21& 1.46& 2.20& 2.21& 0.25& 0.84& 1.10& 1.78& 1.78& 0.34& 0.64& 0.99& 1.54& 1.54\\
 0.4& 0.11& 1.52& 1.75& 2.73& 2.74& 0.16& 1.17& 1.45& 2.25& 2.26& 0.26& 0.78& 1.08& 1.76& 1.76& 0.33& 0.63& 1.00& 1.54& 1.54\\
 0.5& 0.11& 1.36& 1.67& 2.65& 2.66& 0.15& 1.07& 1.35& 2.06& 2.07& 0.25& 0.77& 1.10& 1.76& 1.76& 0.34& 0.65& 1.01& 1.56& 1.56\\
\hline
\end{tabular}
\end{footnotesize}
\end{table}
\end{landscape}

\begin{landscape}
\begin{table}[!h]
\centering
\caption{Relative efficiencies of IPT, PTLE, SLE, and PSLE with respect to LASSO estimator when $n=100$, $r=0.5,$ $p=10$, $\sigma=20$ for different values of $k$.} \label{tb:rel:eff:cor5:p10:sigma20}
\bigskip
\begin{footnotesize}
\begin{tabular}{r|rrrrr|rrrrr|rrrrr|rrrrr}
\hline
 & \multicolumn{5}{c|}{$k=0$} & \multicolumn{5}{c|}{$k=1$} & \multicolumn{5}{c|}{$k=3$}
 & \multicolumn{5}{c}{$k=5$} \\
 $\Delta^2$ & LSE & IPT & PTLE & SLE & PSLE & LSE & IPT & PTLE & SLE & PSLE & 
 LSE & IPT & PTLE & SLE & PSLE& LSE & IPT & PTLE & SLE & PSLE \\
 \hline
 0.0& 0.10& 1.62& 1.91& 3.14& 3.14& 0.09& 1.47& 1.69& 2.75& 2.78& 0.13& 1.15& 1.58& 2.49& 2.50& 0.16& 0.93& 1.30& 2.15& 2.16\\
 0.1& 0.09& 1.59& 1.88& 3.35& 3.36& 0.11& 2.13& 2.45& 3.09& 3.10& 0.13& 1.07& 1.56& 2.38& 2.38& 0.17& 0.89& 1.30& 2.14& 2.15\\
 0.2& 0.11& 1.67& 1.87& 3.06& 3.06& 0.12& 1.62& 1.89& 2.90& 2.90& 0.13& 1.23& 1.74& 2.48& 2.49& 0.16& 0.93& 1.37& 2.19& 2.19\\
 0.3& 0.11& 1.86& 1.97& 3.06& 3.09& 0.09& 1.55& 1.72& 2.80& 2.81& 0.13& 1.12& 1.49& 2.33& 2.33& 0.17& 0.96& 1.38& 2.18& 2.18\\
 0.4& 0.09& 1.83& 2.26& 3.22& 3.24& 0.10& 1.42& 1.85& 2.78& 2.80& 0.13& 1.26& 1.62& 2.41& 2.41& 0.19& 0.98& 1.30& 2.14& 2.15\\
 0.5& 0.10& 1.65& 1.76& 2.94& 2.94& 0.11& 1.51& 1.81& 2.73& 2.74& 0.14& 1.10& 1.46& 2.35& 2.36& 0.19& 0.98& 1.27& 2.20& 2.21\\
\hline
\end{tabular}
\end{footnotesize}
\end{table}
\end{landscape}

\begin{landscape}
\begin{table}[!h]
\centering
\caption{Relative efficiencies of IPT, PTLE, SLE, and PSLE with respect to LASSO estimator when $n=100$, $r=0.5,$ $p=20$, $\sigma=10$ for different values of $k$.} \label{tb:rel:eff:cor5:p20:sigma10}
\bigskip
\begin{footnotesize}
\begin{tabular}{r|rrrrr|rrrrr|rrrrr|rrrrr}
\hline
 & \multicolumn{5}{c|}{$k=0$} & \multicolumn{5}{c|}{$k=1$} & \multicolumn{5}{c|}{$k=3$}
 & \multicolumn{5}{c}{$k=5$} \\
 $\Delta^2$ & LSE & IPT & PTLE & SLE & PSLE & LSE & IPT & PTLE & SLE & PSLE &  
 LSE & IPT & PTLE & SLE & PSLE& LSE & IPT & PTLE & SLE & PSLE \\
 \hline
 0.0& 0.04& 1.50& 2.19& 5.71& 5.77& 0.07& 1.33& 1.93& 3.48& 3.51& 0.10& 0.83& 1.31& 2.25& 2.27& 0.14& 0.57& 1.10& 1.96& 1.96\\
 0.1& 0.05& 1.31& 1.86& 4.35& 4.39& 0.06& 1.31& 2.27& 3.27& 3.31& 0.11& 0.78& 1.26& 2.23& 2.25& 0.15& 0.55& 1.09& 1.99& 2.00\\
 0.2& 0.05& 1.29& 2.04& 4.64& 4.69& 0.07& 1.01& 1.59& 2.92& 2.95& 0.12& 0.76& 1.22& 2.36& 2.37& 0.16& 0.58& 1.07& 1.95& 1.95\\
 0.3& 0.06& 1.57& 2.26& 4.55& 4.59& 0.07& 1.02& 1.51& 2.90& 2.92& 0.11& 0.78& 1.28& 2.25& 2.26& 0.16& 0.56& 1.04& 1.96& 1.96\\
 0.4& 0.06& 1.38& 1.86& 4.10& 4.17& 0.07& 0.99& 1.60& 2.85& 2.88& 0.12& 0.83& 1.26& 2.31& 2.32& 0.16& 0.53& 1.07& 1.92& 1.92\\
 0.5& 0.05& 0.95& 1.62& 3.56& 3.64& 0.07& 1.01& 1.51& 2.68& 2.70& 0.12& 0.75& 1.20& 2.18& 2.19& 0.16& 0.53& 1.05& 1.91& 1.91\\
\hline
\end{tabular}
\end{footnotesize}
\end{table}
\end{landscape}

\begin{landscape}
\begin{table}[!h]
\centering
\caption{Relative efficiencies of IPT, PTLE, SLE, and PSLE with respect to LASSO estimator when $n=100$, $r=0.5,$ $p=20$, $\sigma=20$ for different values of $k$.} \label{tb:rel:eff:cor5:p20:sigma10}
\bigskip
\begin{footnotesize}
\begin{tabular}{r|rrrrr|rrrrr|rrrrr|rrrrr}
\hline
 & \multicolumn{5}{c|}{$k=0$} & \multicolumn{5}{c|}{$k=1$} & \multicolumn{5}{c|}{$k=3$}
 & \multicolumn{5}{c}{$k=5$} \\
 $\Delta^2$ & LSE & IPT & PTLE & SLE & PSLE & LSE & IPT & PTLE & SLE & PSLE & 
 LSE & IPT & PTLE & SLE & PSLE& LSE & IPT & PTLE & SLE & PSLE \\
 \hline
 0.0& 0.05& 1.46& 1.92& 5.45& 5.51& 0.06& 1.32& 1.69& 4.21& 4.25& 0.06& 1.16& 1.72& 3.54& 3.57& 0.08& 0.95& 1.57& 3.04& 3.05\\
 0.1& 0.04& 1.67& 3.53& 7.06& 7.25& 0.05& 1.54& 2.04& 4.36& 4.43& 0.08& 1.23& 1.49& 3.28& 3.30& 0.08& 0.95& 1.60& 3.06& 3.08\\
 0.2& 0.05& 1.39& 1.79& 4.43& 4.50& 0.05& 1.08& 1.78& 4.14& 4.20& 0.06& 1.09& 1.90& 3.67& 3.71& 0.08& 0.88& 1.55& 3.10& 3.13\\
 0.3& 0.05& 1.54& 2.60& 5.84& 5.91& 0.05& 1.33& 1.93& 4.31& 4.35& 0.07& 1.16& 1.59& 3.33& 3.36& 0.08& 0.90& 1.47& 2.88& 2.90\\
 0.4& 0.05& 1.66& 2.63& 5.29& 5.43& 0.06& 1.28& 1.74& 4.21& 4.26& 0.07& 1.23& 1.67& 3.32& 3.38& 0.09& 1.01& 1.44& 3.02& 3.04\\
 0.5& 0.05& 1.25& 2.11& 4.74& 4.88& 0.06& 1.62& 1.87& 4.37& 4.40& 0.07& 1.18& 1.75& 3.48& 3.50& 0.08& 0.89& 1.53& 2.94& 2.97\\
\hline
\end{tabular}
\end{footnotesize}
\end{table}
\end{landscape}

\begin{landscape}
\begin{table}[!h]
\centering
\caption{Relative efficiencies of IPT, PTLE, SLE, and PSLE with respect to LASSO estimator when $n=100$, $r=0.5,$ $p=50$, $\sigma=10$ for different values of $k$.} \label{tb:rel:eff:cor5:p50:sigma10}
\bigskip
\begin{footnotesize}
\begin{tabular}{r|rrrrr|rrrrr|rrrrr|rrrrr}
\hline
 & \multicolumn{5}{c|}{$k=0$} & \multicolumn{5}{c|}{$k=1$} & \multicolumn{5}{c|}{$k=3$}
 & \multicolumn{5}{c}{$k=5$} \\
 $\Delta^2$ & LSE & IPT & PTLE & SLE & PSLE & LSE & IPT & PTLE & SLE & PSLE & 
 LSE & IPT & PTLE & SLE & PSLE& LSE & IPT & PTLE & SLE & PSLE \\
 \hline
 0.0& 0.02& 0.78& 1.58& 7.47& 7.60& 0.02& 0.56& 1.90& 4.49& 4.58& 0.03& 0.59& 1.41& 2.93& 2.98& 0.04& 0.41& 1.28& 2.56& 2.58\\
 0.1& 0.01& 0.62& 2.23& 8.60& 8.95& 0.02& 0.64& 1.58& 4.53& 4.64& 0.03& 0.52& 1.53& 2.97& 3.04& 0.04& 0.40& 1.25& 2.63& 2.65\\
 0.2& 0.01& 0.52& 1.97& 7.00& 7.40& 0.02& 0.57& 1.78& 4.30& 4.39& 0.03& 0.53& 1.43& 3.00& 3.05& 0.04& 0.41& 1.24& 2.59& 2.61\\
 0.3& 0.02& 0.47& 1.63& 5.90& 6.15& 0.02& 0.64& 1.55& 4.03& 4.12& 0.03& 0.41& 1.37& 2.73& 2.79& 0.05& 0.39& 1.21& 2.61& 2.63\\
 0.4& 0.02& 0.55& 1.82& 5.81& 6.04& 0.02& 0.51& 1.42& 3.84& 3.91& 0.03& 0.46& 1.41& 2.74& 2.78& 0.04& 0.36& 1.24& 2.53& 2.55\\
 0.5& 0.02& 0.62& 1.59& 5.31& 5.47& 0.02& 0.56& 1.53& 3.45& 3.50& 0.03& 0.48& 1.37& 2.80& 2.83& 0.05& 0.40& 1.22& 2.54& 2.55\\
\hline
\end{tabular}
\end{footnotesize}
\end{table}
\end{landscape}

\begin{landscape}
\begin{table}[!h]
\centering
\caption{Relative efficiencies of IPT, PTLE, SLE, and PSLE with respect to LASSO estimator when $n=100$, $r=0.5,$ $p=50$, $\sigma=20$ for different values of $k$.} \label{tb:rel:eff:cor5:p50:sigma20}
\bigskip
\begin{footnotesize}
\begin{tabular}{r|rrrrr|rrrrr|rrrrr|rrrrr}
\hline
 & \multicolumn{5}{c|}{$k=0$} & \multicolumn{5}{c|}{$k=1$} & \multicolumn{5}{c|}{$k=3$}
 & \multicolumn{5}{c}{$k=5$} \\
 $\Delta^2$ & LSE & IPT & PTLE & SLE & PSLE & LSE & IPT & PTLE & SLE & PSLE & 
 LSE & IPT & PTLE & SLE & PSLE& LSE & IPT & PTLE & SLE & PSLE \\
 \hline
 0.0& 0.02& 0.60& 2.02& 8.19& 8.76& 0.01& 0.48& 1.63& 5.84& 6.04& 0.02& 0.47& 1.86& 5.08& 5.32& 0.02& 0.56& 1.52& 4.02& 4.13\\
 0.1& 0.01& 0.52& 1.65& 8.29& 8.73& 0.02& 0.56& 1.83& 6.45& 6.66& 0.02& 0.50& 1.56& 4.51& 4.60& 0.02& 0.46& 1.61& 4.36& 4.44\\
 0.2& 0.02& 0.63& 2.02& 7.89& 8.46& 0.02& 0.59& 1.70& 5.92& 6.06& 0.02& 0.54& 1.80& 4.40& 4.51& 0.02& 0.49& 1.81& 4.15& 4.27\\
 0.3& 0.02& 0.56& 1.99& 8.17& 8.57& 0.02& 0.70& 1.64& 6.43& 6.55& 0.02& 0.60& 1.68& 4.75& 4.90& 0.02& 0.56& 1.75& 4.11& 4.20\\
 0.4& 0.01& 0.49& 2.23& 9.43& 9.79& 0.02& 0.60& 1.78& 6.45& 6.58& 0.02& 0.57& 1.57& 4.13& 4.22& 0.02& 0.51& 1.69& 3.79& 3.88\\
 0.5& 0.01& 0.54& 2.04& 9.00& 9.40& 0.01& 0.58& 1.91& 6.46& 6.70& 0.02& 0.59& 1.91& 5.14& 5.29& 0.03& 0.53& 1.41& 3.95& 4.03\\
\hline
\end{tabular}
\end{footnotesize}
\end{table}
\end{landscape}

\begin{landscape}
\begin{table}[!h]
\centering
\caption{Relative efficiencies of IPT, PTLE, SLE, and PSLE with respect to LASSO estimator when $n=100$, $r=0.9,$ $p=10$, $\sigma=10$ for different values of $k$.} \label{tb:rel:eff:cor9:p10:sigma10}
\bigskip
\begin{footnotesize}
\begin{tabular}{r|rrrrr|rrrrr|rrrrr|rrrrr}
\hline
 & \multicolumn{5}{c|}{$k=0$} & \multicolumn{5}{c|}{$k=1$} & \multicolumn{5}{c|}{$k=3$}
 & \multicolumn{5}{c}{$k=5$} \\
 $\Delta^2$ & LSE & IPT & PTLE & SLE & PSLE & LSE & IPT & PTLE & SLE & PSLE & 
 LSE & IPT & PTLE & SLE & PSLE& LSE & IPT & PTLE & SLE & PSLE \\
 \hline
 0.0& 0.09& 1.70& 1.73& 3.09& 3.09& 0.10& 1.47& 1.86& 2.89& 2.90& 0.13& 0.46& 1.14& 1.97& 1.97& 0.18& 0.32& 1.01& 1.56& 1.56\\
 0.1& 0.09& 1.52& 1.83& 3.13& 3.14& 0.09& 1.04& 1.59& 2.65& 2.65& 0.15& 0.45& 1.09& 1.93& 1.93& 0.19& 0.32& 1.00& 1.51& 1.51\\
 0.2& 0.10& 1.54& 1.58& 2.86& 2.87& 0.10& 0.92& 1.59& 2.51& 2.51& 0.14& 0.43& 1.08& 1.94& 1.94& 0.19& 0.32& 1.01& 1.49& 1.49\\
 0.3& 0.08& 1.60& 2.28& 3.27& 3.29& 0.10& 0.99& 1.68& 2.59& 2.60& 0.15& 0.40& 1.10& 1.86& 1.86& 0.20& 0.32& 1.00& 1.46& 1.46\\
 0.4& 0.10& 1.66& 1.57& 2.88& 2.89& 0.11& 0.93& 1.53& 2.54& 2.55& 0.15& 0.42& 1.12& 1.89& 1.89& 0.17& 0.28& 1.00& 1.47& 1.47\\
 0.5& 0.10& 1.44& 1.90& 3.06& 3.07& 0.10& 0.90& 1.51& 2.51& 2.51& 0.14& 0.38& 1.08& 1.82& 1.82& 0.20& 0.32& 1.00& 1.45& 1.45\\
\hline
\end{tabular}
\end{footnotesize}
\end{table}
\end{landscape}

\begin{landscape}
\begin{table}[!h]
\centering
\caption{Relative efficiencies of IPT, PTLE, SLE, and PSLE with respect to LASSO estimator when $n=100$, $r=0.9,$ $p=10$, $\sigma=20$ for different values of $k$.} \label{tb:rel:eff:cor9:p10:sigma20}
\bigskip
\begin{footnotesize}
\begin{tabular}{r|rrrrr|rrrrr|rrrrr|rrrrr}
\hline
 & \multicolumn{5}{c|}{$k=0$} & \multicolumn{5}{c|}{$k=1$} & \multicolumn{5}{c|}{$k=3$}
 & \multicolumn{5}{c}{$k=5$} \\
 $\Delta^2$ & LSE & IPT & PTLE & SLE & PSLE & LSE & IPT & PTLE & SLE & PSLE &  
 LSE & IPT & PTLE & SLE & PSLE& LSE & IPT & PTLE & SLE & PSLE \\
 \hline
 0.0& 0.09& 1.88& 2.14& 3.27& 3.29& 0.09& 1.68& 1.74& 3.00& 3.01& 0.10& 0.99& 1.68& 2.77& 2.78& 0.11& 0.54& 1.30& 2.23& 2.24\\
 0.1& 0.09& 2.43& 2.23& 3.16& 3.17& 0.10& 1.44& 1.60& 2.85& 2.86& 0.10& 0.93& 1.54& 2.62& 2.63& 0.12& 0.55& 1.26& 2.23& 2.23\\
 0.2& 0.10& 1.50& 1.83& 2.98& 3.00& 0.08& 1.28& 1.65& 2.92& 2.92& 0.11& 0.91& 1.47& 2.61& 2.61& 0.13& 0.59& 1.24& 2.29& 2.29\\
 0.3& 0.10& 1.81& 1.93& 3.17& 3.18& 0.10& 1.59& 1.84& 3.08& 3.09& 0.11& 1.00& 1.50& 2.58& 2.59& 0.13& 0.49& 1.17& 2.12& 2.12\\
 0.4& 0.10& 1.85& 2.26& 3.31& 3.32& 0.08& 1.34& 1.92& 3.16& 3.17& 0.10& 0.89& 1.55& 2.54& 2.54& 0.11& 0.51& 1.28& 2.16& 2.16\\
 0.5& 0.10& 1.57& 1.78& 3.19& 3.19& 0.10& 1.37& 1.57& 2.88& 2.90& 0.10& 0.74& 1.43& 2.46& 2.46& 0.13& 0.54& 1.20& 2.16& 2.16\\
\hline
\end{tabular}
\end{footnotesize}
\end{table}
\end{landscape}

\begin{landscape}
\begin{table}[!h]
\centering
\caption{Relative efficiencies of IPT, PTLE, SLE, and PSLE with respect to LASSO estimator when $n=100$, $r=0.9,$ $p=20$, $\sigma=10$ for different values of $k$.} \label{tb:rel:eff:cor9:p20:sigma10}
\bigskip
\begin{footnotesize}
\begin{tabular}{r|rrrrr|rrrrr|rrrrr|rrrrr}
\hline
 & \multicolumn{5}{c|}{$k=0$} & \multicolumn{5}{c|}{$k=1$} & \multicolumn{5}{c|}{$k=3$}
 & \multicolumn{5}{c}{$k=5$} \\
 $\Delta^2$ & LSE & IPT & PTLE & SLE & PSLE & LSE & IPT & PTLE & SLE & PSLE & 
 LSE & IPT & PTLE & SLE & PSLE& LSE & IPT & PTLE & SLE & PSLE \\
 \hline
 0.0& 0.05& 1.26& 1.74& 4.79& 4.82& 0.06& 1.17& 1.55& 3.86& 3.88& 0.06& 0.44& 1.33& 3.08& 3.09& 0.09& 0.28& 1.04& 2.17& 2.18\\
 0.1& 0.05& 1.37& 2.08& 4.99& 5.02& 0.05& 1.19& 1.89& 4.43& 4.47& 0.07& 0.47& 1.33& 3.02& 3.03& 0.10& 0.27& 1.03& 2.04& 2.04\\
 0.2& 0.05& 1.40& 1.80& 4.74& 4.79& 0.05& 0.95& 1.67& 4.24& 4.27& 0.07& 0.41& 1.29& 2.82& 2.83& 0.09& 0.25& 1.03& 1.99& 1.99\\
 0.3& 0.04& 0.93& 1.66& 4.55& 4.59& 0.06& 0.85& 1.49& 3.69& 3.72& 0.07& 0.39& 1.21& 2.71& 2.71& 0.09& 0.23& 1.03& 1.98& 1.98\\
 0.4& 0.04& 1.13& 2.37& 5.19& 5.26& 0.04& 0.80& 1.70& 4.06& 4.10& 0.07& 0.38& 1.22& 2.67& 2.67& 0.10& 0.24& 1.01& 2.01& 2.01\\
 0.5& 0.06& 1.21& 1.53& 3.66& 3.68& 0.05& 0.93& 1.74& 3.80& 3.82& 0.08& 0.36& 1.17& 2.58& 2.59& 0.10& 0.24& 1.02& 1.96& 1.96\\
\hline
\end{tabular}
\end{footnotesize}
\end{table}
\end{landscape}

\begin{landscape}
\begin{table}[!h]
\centering
\caption{Relative efficiencies of IPT, PTLE, SLE, and PSLE with respect to LASSO estimator when $n=100$, $r=0.9,$ $p=20$, $\sigma=20$ for different values of $k$.} \label{tb:rel:eff:cor9:p20:sigma20}
\bigskip
\begin{footnotesize}
\begin{tabular}{r|rrrrr|rrrrr|rrrrr|rrrrr}
\hline
 & \multicolumn{5}{c|}{$k=0$} & \multicolumn{5}{c|}{$k=1$} & \multicolumn{5}{c|}{$k=3$}
 & \multicolumn{5}{c}{$k=5$} \\
 $\Delta^2$ & LSE & IPT & PTLE & SLE & PSLE & LSE & IPT & PTLE & SLE & PSLE &  
 LSE & IPT & PTLE & SLE & PSLE& LSE & IPT & PTLE & SLE & PSLE \\
 \hline
 0.0& 0.04& 1.51& 2.11& 5.76& 5.81& 0.04& 1.12& 1.68& 4.39& 4.42& 0.05& 1.03& 1.80& 4.46& 4.48& 0.05& 0.63& 1.63& 3.69& 3.72\\
 0.1& 0.04& 1.66& 2.35& 5.11& 5.13& 0.05& 1.16& 1.43& 4.26& 4.30& 0.04& 1.01& 2.02& 4.57& 4.62& 0.05& 0.51& 1.44& 3.78& 3.79\\
 0.2& 0.03& 1.12& 2.21& 5.48& 5.71& 0.06& 1.37& 1.65& 4.51& 4.55& 0.05& 0.77& 1.65& 3.91& 3.96& 0.06& 0.61& 1.39& 3.56& 3.58\\
 0.3& 0.04& 1.21& 1.87& 5.11& 5.14& 0.05& 1.22& 1.96& 5.13& 5.17& 0.06& 1.05& 1.94& 4.48& 4.52& 0.07& 0.59& 1.38& 3.34& 3.35\\
 0.4& 0.05& 1.33& 1.97& 5.60& 5.65& 0.05& 1.19& 1.98& 4.86& 4.91& 0.06& 1.02& 1.88& 4.48& 4.50& 0.06& 0.54& 1.37& 3.51& 3.52\\
 0.5& 0.05& 1.23& 2.00& 5.29& 5.35& 0.04& 1.01& 2.06& 4.94& 4.97& 0.05& 0.84& 1.61& 3.81& 3.83& 0.07& 0.53& 1.33& 3.07& 3.08\\
\hline
\end{tabular}
\end{footnotesize}
\end{table}
\end{landscape}

\begin{landscape}
\begin{table}[!h]
\centering
\caption{Relative efficiencies of IPT, PTLE, SLE, and PSLE with respect to LASSO estimator when $n=100$, $r=0.9,$ $p=50$, $\sigma=10$ for different values of $k$.} \label{tb:rel:eff:cor9:p50:sigma10}
\bigskip
\begin{footnotesize}
\begin{tabular}{r|rrrrr|rrrrr|rrrrr|rrrrr}
\hline
 & \multicolumn{5}{c|}{$k=0$} & \multicolumn{5}{c|}{$k=1$} & \multicolumn{5}{c|}{$k=3$}
 & \multicolumn{5}{c}{$k=5$} \\
 $\Delta^2$ & LSE & IPT & PTLE & SLE & PSLE & LSE & IPT & PTLE & SLE & PSLE &  LSE & IPT & PTLE & SLE & PSLE& LSE & IPT & PTLE & SLE & PSLE \\
 \hline
 0.0& 0.01& 0.50& 2.72& 12.93& 13.75& 0.02& 0.57& 1.93& 6.34& 6.62& 0.02& 0.29& 1.50& 4.44& 4.49& 0.03& 0.20& 1.16& 3.47& 3.48\\
 0.1& 0.02& 0.54& 2.03&  9.09&  9.46& 0.01& 0.42& 1.77& 6.33& 6.58& 0.02& 0.28& 1.39& 4.63& 4.68& 0.03& 0.17& 1.13& 3.26& 3.26\\
 0.2& 0.01& 0.55& 2.17&  7.98&  8.39& 0.02& 0.42& 1.84& 6.48& 6.80& 0.02& 0.28& 1.40& 4.62& 4.67& 0.02& 0.15& 1.12& 3.22& 3.22\\
 0.3& 0.01& 0.52& 1.94&  7.59&  7.84& 0.01& 0.53& 2.31& 7.46& 7.83& 0.02& 0.28& 1.46& 4.80& 4.87& 0.02& 0.14& 1.17& 3.37& 3.37\\
 0.4& 0.02& 0.49& 1.93&  5.46&  5.53& 0.02& 0.50& 2.08& 7.17& 7.41& 0.02& 0.25& 1.31& 3.98& 4.00& 0.03& 0.16& 1.14& 3.17& 3.17\\
 0.5& 0.01& 0.52& 2.05&  7.33&  7.66& 0.02& 0.38& 1.52& 5.92& 6.12& 0.02& 0.25& 1.51& 4.37& 4.41& 0.03& 0.15& 1.11& 3.17& 3.17\\
\hline
\end{tabular}
\end{footnotesize}
\end{table}
\end{landscape}

\begin{landscape}
\begin{table}[!h]
\centering
\caption{Relative efficiencies of IPT, PTLE, SLE, and PSLE with respect to LASSO estimator when $n=100$, $r=0.9,$ $p=50$, $\sigma=20$ for different values of $k$.} \label{tb:rel:eff:cor9:p50:sigma20}
\bigskip
\begin{footnotesize}
\begin{tabular}{r|rrrrr|rrrrr|rrrrr|rrrrr}
\hline
 & \multicolumn{5}{c|}{$k=0$} & \multicolumn{5}{c|}{$k=1$} & \multicolumn{5}{c|}{$k=3$}
 & \multicolumn{5}{c}{$k=5$} \\
 $\Delta^2$ & LSE & IPT & PTLE & SLE & PSLE & LSE & IPT & PTLE & SLE & PSLE & 
 LSE & IPT & PTLE & SLE & PSLE& LSE & IPT & PTLE & SLE & PSLE \\
 \hline
 0.0& 0.01& 0.55& 2.34& 10.73& 11.28& 0.01& 0.53& 2.17& 9.25& 9.67& 0.02& 0.51& 1.99& 7.74& 7.99& 0.02& 0.37& 1.72& 6.29& 6.52\\
 0.1& 0.01& 0.45& 2.35&  9.83& 10.31& 0.02& 0.53& 1.61& 6.96& 7.12& 0.02& 0.46& 1.76& 7.18& 7.32& 0.02& 0.33& 1.52& 5.64& 5.75\\
 0.2& 0.01& 0.44& 1.92&  9.03&  9.34& 0.01& 0.53& 2.16& 9.44& 9.98& 0.02& 0.42& 1.71& 5.85& 6.03& 0.02& 0.34& 1.61& 5.61& 5.74\\
 0.3& 0.01& 0.49& 2.20&  9.35&  9.84& 0.01& 0.48& 1.87& 7.93& 8.36& 0.02& 0.46& 1.90& 6.61& 6.77& 0.02& 0.32& 1.62& 6.08& 6.19\\
 0.4& 0.01& 0.55& 2.00&  8.62&  8.94& 0.02& 0.50& 1.78& 7.63& 7.82& 0.02& 0.56& 1.84& 7.30& 7.49& 0.02& 0.39& 1.58& 5.74& 5.85\\
 0.5& 0.01& 0.47& 1.93&  8.90&  9.33& 0.01& 0.52& 2.06& 7.89& 8.11& 0.02& 0.43& 1.78& 6.74& 6.96& 0.02& 0.34& 1.63& 5.81& 5.92\\
\hline
\end{tabular}
\end{footnotesize}
\end{table}
\end{landscape}

\bibliographystyle{apa}
\bibliography{ShrinkageReferences}

\begin{thebibliography}{}

\bibitem[\protect\astroncite{Bancroft}{1944}]{bancroft:1944}
Bancroft, T.~A. (1944).
\newblock On biases in estimation due to the use of preliminary tests of
  significances.
\newblock {\em Annals of Mathematical Statistics}, 15.

\bibitem[\protect\astroncite{Donoho and Johnstone}{1994}]{donoho:1994:ideal}
Donoho, D.~L. and Johnstone, J.~M. (1994).
\newblock Ideal spatial adaptation by wavelet shrinkage.
\newblock {\em Biometrika}, 81(3):425--455.

\bibitem[\protect\astroncite{Efron et~al.}{2004}]{efron:etal:2004}
Efron, B., Hastie, T., Johnstone, I., and Tibshirani, R. (2004).
\newblock Least angle regression.
\newblock {\em Annals of Statistics}, 32:407--499.

\bibitem[\protect\astroncite{Fan and Li}{2001}]{fan:li:2001}
Fan, J. and Li, R. (2001).
\newblock Variable selection via nonconcave penalized likelihood and its oracle
  properties.
\newblock {\em Journal of the American Statistical Association},
  96(456):1348--1360.

\bibitem[\protect\astroncite{Faraway}{2002}]{faraway:2002}
Faraway, J.~J. (2002).
\newblock {\em Practical Regression and Anova using R}.
\newblock Website.

\bibitem[\protect\astroncite{Frank and Friedman}{1993}]{frank:friedman:1993}
Frank, I.~E. and Friedman, J.~H. (1993).
\newblock A statistical view of some chemometrics regression tools.
\newblock {\em Technometrics}, 35:109--148.

\bibitem[\protect\astroncite{Friedman et~al.}{2014}]{glmnet:package}
Friedman, J., Hastie, T., Simon, N., and Tibshirani, R. (2014).
\newblock {\em glmnet: fit a GLM with lasso or elasticnet regularization}.
\newblock R package version 1.9-8 --- For new features, see the `Changelog'
  file (in the package source).

\bibitem[\protect\astroncite{Friedman
  et~al.}{2010}]{friedman:hastie:tibshirani:2010}
Friedman, J., Hastie, T., and Tibshirani, R. (2010).
\newblock Regularization paths for generalized linear models via coordinate
  descent.
\newblock {\em Journal of Statistical Software}, 33(1):1--22.

\bibitem[\protect\astroncite{James and Stein}{1961}]{james:stein:1961}
James, W. and Stein, C. (1961).
\newblock Estimation with quadratic loss.
\newblock In {\em Proceedings of the Fourth Berkeley Symposium on Mathematical
  Statistics and Probability, Volume 1: Contributions to the Theory of
  Statistics}, pages 361--379, Berkeley, Calif. University of California Press.

\bibitem[\protect\astroncite{Judge and Bock}{1978}]{judge:bock:1978}
Judge, G.~G. and Bock, M.~E. (1978).
\newblock {\em The Statistical Implications of Pre-test and Stein-rule
  Estimators in Econometrics.}
\newblock North Holland, Amsterdam.

\bibitem[\protect\astroncite{Longley}{1967}]{longley:1967}
Longley, J.~W. (1967).
\newblock An appraisal of least squares programs for the electronic computer
  from the viewpoint of the user.
\newblock {\em Journal of the American Statistical Association}, 62:819--841.

\bibitem[\protect\astroncite{{R Core Team}}{2014}]{r:2014}
{R Core Team} (2014).
\newblock {\em R: A Language and Environment for Statistical Computing}.
\newblock R Foundation for Statistical Computing, Vienna, Austria.

\bibitem[\protect\astroncite{Raheem et~al.}{2012}]{raheem:ahmed:doksum:2012}
Raheem, E., Ahmed, S.~E., and Doksum, K.~A. (2012).
\newblock Absolute penalty and shrinkage estimation in partially linear models.
\newblock {\em Computational Statistics \& Data Analysis}, 56(4):874--891.

\bibitem[\protect\astroncite{Saleh}{2006}]{saleh:2006}
Saleh, A. K. M.~E. (2006).
\newblock {\em Theory of Preliminary Test and Stein-Type Estimation with
  Applications}.
\newblock Wiley.

\bibitem[\protect\astroncite{Saleh and Raheem}{2015}]{saleh:raheem:2015}
Saleh, A. K. M.~E. and Raheem, E. (2015).
\newblock Penalty, shrinkage, and preliminary test estimators under full model
  hypothesis.
\newblock Submitted.

\bibitem[\protect\astroncite{Stein}{1956}]{stein:1956}
Stein, C. (1956).
\newblock The admissibility of hotelling's $t\sp 2$-test.
\newblock {\em Mathematical Statistics}, 27:616--623.

\bibitem[\protect\astroncite{Tibshirani}{1996}]{tibshirani:1996}
Tibshirani, R. (1996).
\newblock Regression shrinkage and selection via the lasso.
\newblock {\em Journal of the Royal Statistical Society: Series B},
  58:267--288.

\end{thebibliography}
\renewcommand{\leftmark}{References} 

 \end{document}